\documentclass[12pt]{article}
\usepackage{bbm}
\usepackage{latexsym}
\usepackage{mathrsfs}
\usepackage{graphicx}
\usepackage{subfigure}
\usepackage{amssymb}
\usepackage{amsmath}
\usepackage{amsthm}
\usepackage{color}
\usepackage{cite}
\usepackage{float}
\usepackage{indentfirst}
\usepackage{anysize}\marginsize{45mm}{45mm}{40mm}{50mm}
\usepackage{caption}
\usepackage{float}
\usepackage{multirow}

\usepackage{enumerate}

\usepackage{latexsym, bm}
\newtheoremstyle{lemma}{\topsep}{\topsep}%
     {}
     {}
     {\bfseries}
     {}
     {0.1em}
     {\thmname{#1}\thmnumber{ #2}\thmnote{ #3}}
\theoremstyle{lemma}  

\newtheorem{theorem}{Theorem}     
\newtheorem{lemma}[theorem]{Lemma}
\newtheorem{corollary}[theorem]{Corollary}

\newtheorem{definition}{Definition}

\numberwithin{equation}{section}

\title{ Paired many-to-many 2-disjoint path cover of balanced hypercubes with faulty edges}

\author{ Huazhong L\"{u}\thanks{
\emph{E-mail address}: lvhz08@lzu.edu.cn} \\
{\small School of Mathematics Science, } \\
{\small University of Electronic Science and Technology of China,}\\
{\small Chengdu 610054, P.R. China}\\
}
\date{}
\begin{document}

\maketitle
\begin{abstract}

As a variant of the well-known hypercube, the balanced hypercube $BH_n$ was proposed as a novel interconnection network topology for parallel computing. It is known that $BH_n$ is bipartite. Assume that $S=\{s_1,s_2\}$ and $T=\{t_1,t_2\}$ are any two sets of two vertices in different partite sets of $BH_n$ ($n\geq1$). It has been proved that there exist two vertex-disjoint $s_1,t_1$-path and $s_2,t_2$-path of $BH_n$ covering all vertices of it. In this paper, we prove that there always exist two vertex-disjoint $s_1,t_1$-path and $s_2,t_2$-path covering all vertices of $BH_n$ with at most $2n-3$ faulty edges. The upper bound $2n-3$ of edge faults tolerated is optimal.

\vskip 0.1 in

\end{abstract}

\section{Introduction}
The interconnection network (network for short) plays a crucial role in massively parallel systems \cite{Leighton}. It is impossible to design a network which is optimum in all aspects of performance, accordingly, many networks have been proposed. Linear arrays and rings are two fundamental networks. Since some parallel applications such as those in image and signal processing are originally designated on an array architecture, it is important to have effective path embedding in a network.

In path embedding problems, to find parallel paths among vertices in networks is one of most central issues concerned with efficient data transmission \cite{Leighton}. Parallel paths in networks are usually studied with regard to disjoint paths in graphs. Since algorithms designed on linear arrays or rings can be efficiently simulated in a topology containing Hamiltonian paths or cycles, Hamiltonian path and cycle embedding property of graphs have been widely studied \cite{Chen,Cai,Cai2,Dybizbanski,Gould,Wang,Xu,Yan}.

In disjoint path cover problems, the many-to-many disjoint path cover problem is the most generalized one\cite{Park2}. Assume that $S=\{s_1, s_2,\cdots, s_k\}$ and $T=\{t_1, t_2,\cdots, t_k\}$ are two sets of $k$ sources and $k$ sinks in a graph $G$, respectively, the {\em many-to-many $k$-disjoint path cover} ($k$-DPC for short) problem is to determine whether there exist $k$ disjoint paths $P_1,P_2,\cdots,P_k$ in $G$ such that $P_i$ joins $s_i$ to $t_i$ for each $i\in\{1,2,\cdots,k\}$ and $V(P_1)\cup\cdots \cup V(P_k)=V(G)$. Moreover, the DPC problem has a close relationship with Hamiltonian path problem in graphs. In fact, a one-to-one DPC of a network is indeed a Hamiltonian path between any two vertices.

Failure is inevitable when a massive system is put in use, so it is of great practical importance to consider the fault-tolerant capacity of a network. Hamiltonicity and $k$-DPC problems of various networks with faulty elements were investigated in literature, for example, $k$-ary $n$-cubes \cite{Chen,Wang2}, recursive circulants \cite{Tsai2,Kim}, hypercubes \cite{Jo,Tsai,Wang} and hypercube-like graphs \cite{Park,Dong}.

The balanced hypercube, proposed by Wu and Huang \cite{Wu}, is one of the most popular networks. It has many excellent topological properties, such as high symmetry, low-latency, regularity, strong connectivity, etc. The special property of the balanced hypercube is that each processor has a backup processor that shares the same neighborhood. Thus tasks running on a faulty processor can be shifted to its backup processor \cite{Wu}. With such novel properties above, different aspects of the balanced hypercube are studied extensively,
including Hamiltonian embedding issues \cite{Cheng,Hao,Li,Lu4,Xu,Yang,Zhou}, connectivity issues \cite{Lu2,Yang3},
matching preclusion and extendability \cite{Lu,Lu3}, and symmetric issues \cite{Zhou2,Zhou3} and some other topics \cite{Huang,Yang2}. In this paper, we will consider the problem of paired 2-DPC of the balanced hypercube with faulty edges.

The rest of this paper is organized as follows. In Section 2, some definitions and notations are presented. The main result of this paper is shown in Section 3. Conclusions are given in Section 4.

\section{Definitions and preliminaries}

Throughout this paper, a network is represented by a simple undirected graph, where vertices represent processors and edges represent links between processors. Let $G=(V(G),E(G))$ be a graph, where $V(G)$ and $E(G)$ are its vertex-set and edge-set, respectively. The number of vertices of $G$ is denoted by $|V(G)|$. The set of vertices adjacent to $v$ is called {\em neighborhood} of $v$, denoted by $N_G(v)$. We will use $N(v)$ to replace $N_G(v)$ when the context is clear. A path $P$ in $G$ is a sequence of distinct vertices so that there is an edge joining consecutive vertices, and the length of $P$ is the number of edges, denoted by $l(P)$. For simplicity, a path $P=\langle x_0,x_1,\cdots, x_k\rangle$ can also be denoted by $\langle x_0,P,x_k\rangle$. A $u,v$-path is a path whose end vertices are $u$ and $v$. If a path $C=\langle x_0,x_1,\cdots, x_k\rangle$ is such that $k\geq3$, $x_0=x_k$, then $C$ is said to be a {\em cycle}, and the length of $C$ is the number of edges. The {\em distance} between two vertices $u$ and $v$, denoted by $d(u, v)$, is the length of a shortest path of $G$ joining $u$ and $v$. A path (resp. cycle) containing all vertices of a graph $G$ is called a {\em Hamiltonian path} (resp. {\em cycle}). A bipartite graph $G$ is {\em bipanconnected} if, for two arbitrary nodes $u$ and $v$ of $G$ with distance $d(u,v)$, there exists a path of length $l$ between $u$ and $v$ for every integer $l$ with $d(u,v)\leq l\leq |V(G)|-1$ and $l\equiv d(u,v)($mod $2)$. For other standard graph notations not defined here please refer to \cite{Bondy}.
\vskip 0.05 in

The definitions of the balanced hypercube are given as follows.
\vskip 0.0 in

\begin{definition}{\bf .}\label{def1}\cite{Wu} An $n$-dimension balanced hypercube $BH_{n}$ contains
$4^{n}$ vertices $(a_{0},$ $\ldots,a_{i-1},$
$a_{i},a_{i+1},\ldots,a_{n-1})$, where $a_{i}\in\{0,1,2,3\},$ $0\leq
i\leq n-1$. Any vertex $v=(a_{0},\ldots,a_{i-1},$
$a_{i},a_{i+1},\ldots,a_{n-1})$ in $BH_{n}$ has the following $2n$ neighbors:

\begin{enumerate}
\item $((a_{0}+1)$ mod $
4,a_{1},\ldots,a_{i-1},a_{i},a_{i+1},\ldots,a_{n-1})$,\\
      $((a_{0}-1)$ mod $ 4,a_{1},\ldots,a_{i-1},a_{i},a_{i+1},\ldots,a_{n-1})$, and
\item $((a_{0}+1)$ mod $ 4,a_{1},\ldots,a_{i-1},(a_{i}+(-1)^{a_{0}})$ mod $
4,a_{i+1},\ldots,a_{n-1})$,\\
      $((a_{0}-1)$ mod $ 4,a_{1},\ldots,a_{i-1},(a_{i}+(-1)^{a_{0}})$ mod $
      4,a_{i+1},\ldots,a_{n-1})$.
\end{enumerate}
\end{definition}

The first coordinate $a_{0}$ of the vertex
$(a_{0},\ldots,a_{i},\ldots,a_{n-1})$ in $BH_{n}$ is defined as {\em inner index}, and
other coordinates $a_{i}$ $(1\leq i\leq n-1)$ {\em outer index}.

\vskip 0.0 in

The recursive structure of the balanced hypercube is presented in the following definition.

\begin{definition}{\bf .}\label{def2}\cite{Wu}
\begin{enumerate}
\item $BH_{1}$ is a $4$-cycle, whose vertices are labelled
by $0,1,2,3$ clockwise.
\item $BH_{k+1}$ is constructed from $4$ $BH_{k}$s, which
are labelled by $BH^{0}_{k}$, $BH^{1}_{k}$, $BH^{2}_{k}$,
$BH^{3}_{k}$. For any vertex in $BH_{k}^{i}(0\leq i\leq 3)$, its new labelling in $BH_{k+1}$ is $(a_{0},a_{1},\ldots,a_{k-1},i)$, and it has two new neighbors:
\begin{enumerate}
\item[a)] $BH^{i+1}_{k}:((a_{0}+1)$mod $4,a_{1},\ldots,a_{k-1},(i+1)$mod $4)$ and

$((a_{0}-1)$mod $4,a_{1},\ldots,a_{k-1},(i+1)$mod $4)$ if $a_{0}$ is even.

\item[b)] $BH^{i-1}_{k}:((a_{0}+1)$mod $4,a_{1},\ldots,a_{k-1},(i-1)$mod $4)$ and

$((a_{0}-1)$mod $4,a_{1},\ldots,a_{k-1},(i-1)$mod $4)$ if $a_{0}$ is odd.
\end{enumerate}

\end{enumerate}
\end{definition}

\begin{figure}
\centering
\includegraphics[width=150mm]{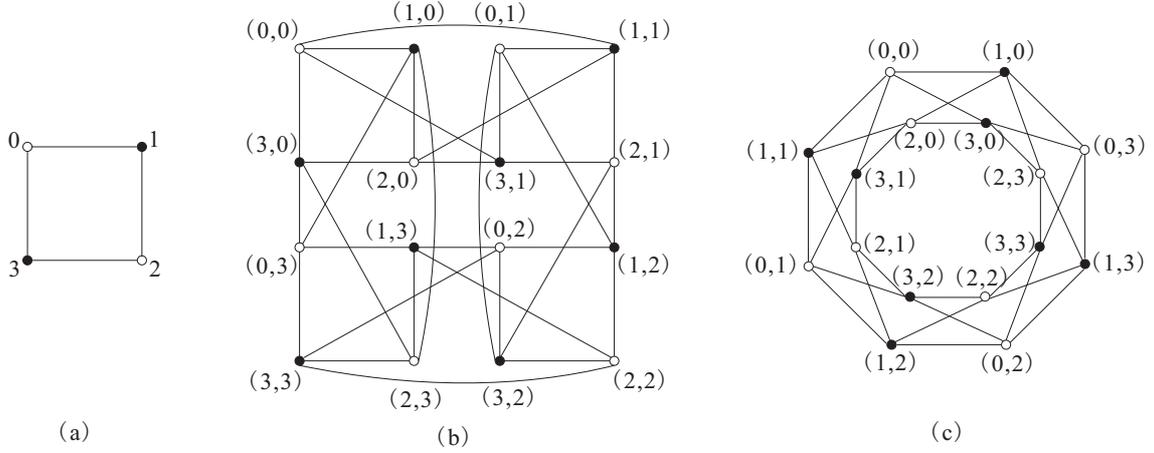}
\caption{$BH_{1}$ and $BH_{2}$.} \label{g1}
\end{figure}

$BH_{1}$ is shown in Fig. \ref{g1} (a). One layout of $BH_{2}$ is shown in Fig. \ref{g1} (b) and the other layout of $BH_{2}$ is shown in Fig. \ref{g1} (c), which reveals a ring-like structure of $BH_2$. Obviously, $BH_{2}$ can be also regarded as joining diagonal vertices of eight twisted 4-cycles end-to-end.

The following basic properties of the balanced hypercube will be applied in the main result of this paper.

\begin{lemma}\label{bipartite}\cite{Wu}{\bf.}
$BH_{n}$ is bipartite.
\end{lemma}

By the above lemma, we give a bipartition $V_0$ and $V_1$ of $BH_{n}$, where $V_0=\{(a_0,\cdots,a_{n-1})|(a_0,\cdots,a_{n-1})\in V(BH_n)$ and $a_0$ is even$\}$ and $V_1=\{(a_0,\cdots,a_{n-1})|$ $(a_0,\cdots,a_{n-1})\in V(BH_n)$ and $a_0$ is odd$\}$.

\begin{lemma}\label{transitive}\cite{Wu,Zhou}{\bf.}
$BH_{n}$ is vertex-transitive and edge-transitive.
\end{lemma}

\begin{lemma}\label{neighbor}\cite{Wu}{\bf.}
Vertices $u=(a_{0},a_{1},\ldots,a_{n-1})$ and
$v=((a_{0}+2)$ mod 4, $a_{1},\ldots,a_{n-1})$ in $BH_{n}$ have the same neighborhood.
\end{lemma}

For convenience, let $p(u)$ be the vertex having the same neighborhood of $u$. It is obvious that $u$ and $p(u)$ differ only from the inner index.

Assume that $u$ is a neighbor of $v$ in $BH_n$. If $u$ and $v$ differ only from the inner index, then $uv$ is called a $0$-{\em dimension edge}, and $u$ and $v$ are mutually called 0-dimension neighbors. Similarly, if $u$ and $v$ differ from $j$-th outer index ($1\leq j\leq n-1$), $uv$ is called a $j$-{\em dimension edge}, and $u$ and $v$ are mutually called $j$-dimension neighbors. The set of all $k$-dimension edges of $BH_n$ is denoted by $E_k$ for each $k\in\{0,\cdots,n-1\}$, and the subgraph of $BH_{n}$ obtained by deleting $E_{n-1}$ is written by $B^{i}$, where $i$ $0\leq i \leq 3$. Obviously, each of $B^{i}$ is isomorphic to $BH_{n-1}$. Let $u_i,v_i,w_i\in V_0$ (resp. $a_i,b_i,c_i\in V_1$) be vertices in $B^i$. For convenience, let $E_{i,i+1}$ be the edge set containing all edges between $B^{i}$ and $B^{i+1}$ ($0\leq i\leq 3$), where ``+'' is under modulo four. For any vertex $v$ of $BH_n$, let $e(v)$ be the set of edges incident to $v$. In particular, the two $k$-dimension edges incident to $v$ is denoted by $e_k(v)$, where $0\leq k\leq n-1$. Let $F$ be a set of edges in $BH_n$, we denote $F^i=F\cap E(B^i)$.

\begin{figure}
\centering
\includegraphics[width=150mm]{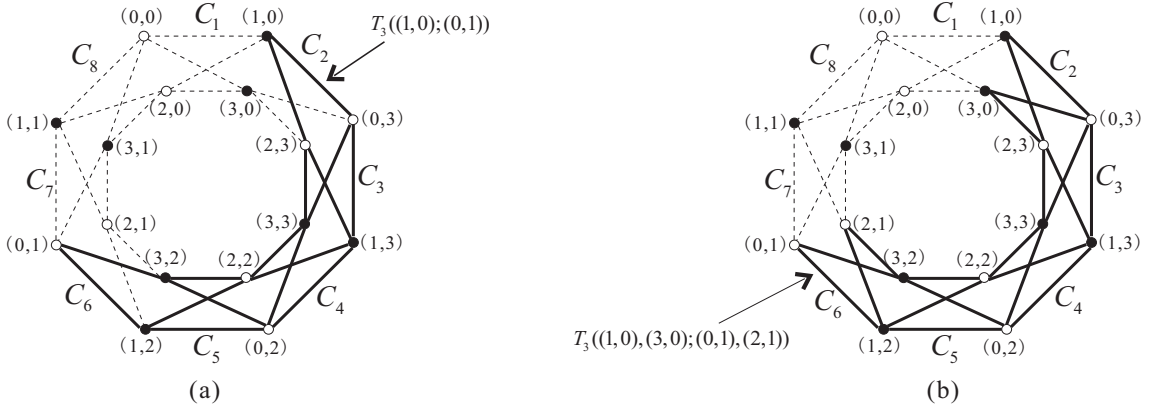}
\caption{$T((1,0);(0,1))$ and $T((1,0),(3,0);(0,1),(3,1))$.} \label{g_Tenon}
\end{figure}

Let $P$ and $Q$ be two 2-paths with central vertices $u$ and $v$, respectively. A {\em tenon chain} $T_m(x;y)$ from $u$ to $v$ is defined to be an $m$ ($m\geq1$) twisted 4-cycle chain with $P$ and $Q$ joining to its two ends, respectively. Additionally, let $P'$ and $Q'$ be two 2-paths with central vertices $x$ and $y$, respectively. $P'$ and $Q'$ are joined to two ends of $T_m(x;y)$ the same way as $P$ and $Q$ do, we denote the graph obtained above by $T_m(u,x;v,y)$. In other words, $T_m(u,x;v,y)$ is an $m+2$ ($m\geq1$) twisted 4-cycles chain with $u$ and $x$ being degree 2 vertices at one end and $v$ and $y$ being degree 2 vertices at the other end. By above, if $1\leq m\leq 6$, $T_m(u;v)$ and $T_m(u,x;v,y)$ are both subgraphs of $BH_2$. For convenience, we refer $T_m(u;v)$ and $T_m(u,x;v,y)$ ($1\leq m\leq6$) to the subgraph of $BH_2$ (ring-like layout) from $u$ to $v$ clockwise. $T_3((1,0),(0,1))$ and $T_3((1,0),(3,0);(0,1),(2,1))$ are illustrated as heavy lines in Figs. \ref{g_Tenon} (a) and \ref{g_Tenon} (b), respectively. Note that if $u$ and $v$ are in different partite sets of $BH_2$ then $m$ is odd, otherwise, $m$ is even.

\section{Paired two-disjoint path cover of the balanced hypercube}

Firstly, we will give some statements, which will be used later.

\begin{lemma}{\rm \cite{Yang2}}{\bf.}\label{common-nb}
Let $u$ be an arbitrary vertex of $BH_n$ for $n\geq1$. Then, for an arbitrary vertex $v$ of $BH_n$, either $u$ and $v$ have 0, 2, or $2n$ common neighbors. Furthermore, there is exactly one vertex $w$ such that $u$ and $w$ have $2n$ common neighbors.
\end{lemma}

\begin{lemma}\label{bipan}
{\rm \cite{Yang}}{\bf.} The balanced hypercube $BH_{n}$ is
bipanconnected for all $n\geq 1$.
\end{lemma}


\begin{lemma}\label{crossing-edges}{\rm \cite{Yang3}}{\bf.} Assume that $n\geq2$. There exist $4^{n-1}$ edges between $B^i$ and $B^{i+1}$ for each $0\leq i\leq3$.
\end{lemma}

\begin{lemma}\cite{Xu}\label{8cycle}{\bf .} Let $uv$ be an edge of $BH_{n}$. Then $uv$ is contained in a cycle $C$ of length 8 in $BH_{n}$ such that $|E(C)\cap E(B^i)|=1$ for each $i=0,1,2,3$.
\end{lemma}

\begin{lemma}{\rm \cite{Cheng}}\label{2paths}{\bf .}  Let $u,x\in V_{0}$ and $v,y\in V_{1}$. Then there exist two vertex-disjoint paths $P$ and $Q$ such that: (1) $P$ connects $u$ to $v$, (2) $Q$ connects $x$ to $y$, (3) $V(P)\cup V(Q)=V(BH_{n})$.
\end{lemma}

\begin{lemma} {\rm \cite{Zhou}}\label{HP-fault}{\bf .} Let $F$ be a set of faulty edges of $BH_n$ with $|F|\leq 2n-2$ for $n\geq 2$ and let $x$ and $y$ be two vertices in different partite sets of $BH_n$. Then there exists a Hamiltonian path of $BH_n-F$ from $x$ to $y$.
\end{lemma}

\begin{lemma} \label{T_xy}{\bf .} Given $T_m(x,y)$ with $m$ being odd. If $f$ is an arbitrary edge of $T_m(x,y)$, then there exists a Hamiltonian path of $T_m(x,y)-f$ from $x$ to $y$.
\end{lemma}

\noindent{\bf Proof.} Since $m$ is odd, $x$ and $y$ are in different partite sets. Either $f$ is an edge incident to $x$ or $y$, or $f$ is an edge of any twisted 4-cycle, it is easy to obtain a Hamiltonian path of $T_m(x,y)$ avoiding $f$. The lemma holds. \qed

It follows from Lemma \ref{T_xy} that there exists a Hamiltonian path of $T_m(x,y)$ from $x$ to $y$ when at most one edge fault occurs, so we also use $T_m(x,y)$ to denote a fault-free Hamiltonian path of $T_m(x,y)$ from $x$ to $y$ when there is no ambiguity.


\begin{lemma} \label{T_uvxy}{\bf .} Given $T_m(u,x;v,y)$ with $m$ being odd. Let $e$ and $f$ be two edges of $T_m(u,x;v,y)$ such that $e$ and $f$ are not contained in the same twisted 4-cycle, then there exist vertex-disjoint $u,v$-path and $x,y$-path of $T_m(u,x;v,y)-\{e,f\}$ that cover all vertices of it.
\end{lemma}

\noindent{\bf Proof.} Since $m$ is odd, $u$ and $x$ are in one partite set, and $v$ and $y$ are in the other partite set of $T_m(u,x;v,y)$. To obtain the desired $u,v$-path and $x,y$-path, one has to go through all twisted 4-cycles of $T_m(u,x;v,y)$ and never go back. Accordingly, $u,v$-path and $x,y$-path contain the same number of vertices. Fault-free $u,v$-path and $x,y$-path of $T_m(u,x;v,y)-\{e,f\}$ can be constructed according to the following two rules:

(1) If $e$ (or $f$) is incident to one of $u,x,v$ and $y$, say $u$, we then choose the other edge incident to $u$ in $u,v$-path.

(2) If $e=ab$ (or $f=ab$) is contained in a twisted 4-cycle $C=\langle a,b,c,d,a\rangle$, then $ad$ (resp. $bc$) must be contained in exact one of $u,v$-path and $x,y$-path.

Hence, the lemma holds. \qed

\begin{lemma}\label{base-main-e-f}{\bf.} Let $\{s_1,s_2\}$ and $\{t_1,t_2\}$ be two sets of vertices in different partite sets of $BH_{2}$ and let $F=\{e,f\}$ be a set of edges of $BH_2$ with $e\in E_0$ and $f\in E_1$. Then there exist vertex-disjoint $s_1,t_1$-path and $s_2,t_2$-path of $BH_{2}-F$ that cover all vertices of it unless there exists a common neighbor of $s_1$ and $s_2$ (or $t_1$ and $t_2$), say $x$, such that $F=e(x)\setminus\{s_1x,s_2x\}$ (or $F=e(x)\setminus\{t_1x,t_2x\}$).
\end{lemma}
\noindent{\bf Proof.} Suppose without loss of generality that $x$ is a common neighbor of $s_1$ and $s_2$, if $F=e(x)\setminus\{s_1x,s_2x\}$, that is, $\{s_1x,s_2x\}\cap F=\emptyset$, which yields a 2-path starting from $s_1$ to $s_2$. Accordingly, it is impossible to obtain vertex-disjoint $s_1,t_1$-path and $s_2,t_2$-path that cover all vertices of $BH_2$. If $d(s_1,s_2)=2$, $F\neq e(x)\setminus\{s_1x,s_2x\}$ is a necessary condition to guarantee that there exist vertex-disjoint $s_1,t_1$-path and $s_2,t_2$-path of $BH_{2}-F$.

On the other hand, noting $e\in E_0$ and $f\in E_1$, each twisted 4-cycle of $BH_{2}$ (ring-like layout) contains at most one of them. By vertex-transitivity of $BH_2$, we may assume that $s_1=(0,0)$. According to all possible relative positions of $s_1,s_2,t_1$ and $t_2$ in $BH_2$, there are 16 essential different distributions to be considered. In each case, we have verified that there always exist vertex-disjoint $s_1,t_1$-path and $s_2,t_2$-path of $BH_{2}-F$ that covers all vertices of $BH_{2}$ (by making use of Lemmas \ref{T_xy} and \ref{T_uvxy} to reduce the number of cases to be cosidered). Since the proof is tedious and rather long, we only list all different distributions of $s_1,s_2,t_1$ and $t_2$ in $BH_2$ in the following:


%

(1) $s_2=(2,0), t_1=(1,0), t_2=(3,0)$; 

(2) $s_2=(2,0), t_1=(1,0), t_2=(3,3)$ (or $s_2=(2,1), t_1=(1,0), t_2=(3,0)$ or $s_2=(2,0), t_1=(3,3), t_2=(3,0)$ or $s_2=(2,3), t_1=(3,3), t_2=(1,3)$); 

(3) $s_2=(2,0), t_1=(1,0), t_2=(3,0)$ (or $s_2=(2,2), t_1=(1,0), t_2=(3,0)$ or $s_2=(2,0), t_1=(3,3), t_2=(3,1)$ or $s_2=(2,2), t_1=(3,3), t_2=(1,3)$); 

(4) $s_2=(2,0), t_1=(1,0), t_2=(3,1)$ (or $s_2=(2,3), t_1=(1,0), t_2=(3,0)$);   

(5) $s_2=(2,3), t_1=(1,0), t_2=(3,3)$ (or $s_2=(2,1), t_1=(1,0), t_2=(3,1)$); 

(6) $s_2=(2,3), t_1=(1,0), t_2=(3,2)$ (or $s_2=(2,2), t_1=(1,0), t_2=(3,1)$ or $s_2=(2,2), t_1=(3,3), t_2=(3,2)$ or $s_2=(2,1), t_1=(3,3), t_2=(3,1)$);  

(7) $s_2=(2,3), t_1=(1,0), t_2=(3,1)$ (or $s_2=(2,3), t_1=(3,3), t_2=(3,0)$); 

(8) $s_2=(2,2), t_1=(1,0), t_2=(3,3)$ (or $s_2=(2,1), t_1=(1,0), t_2=(3,2)$); 

(9) $s_2=(2,2), t_1=(1,0), t_2=(3,2)$; 

(10) $s_2=(2,1), t_1=(1,0), t_2=(3,3)$ (or $s_2=(2,1), t_1=(3,3), t_2=(3,2)$); 

(11) $s_2=(2,0), t_1=(1,3), t_2=(3,3)$;  

(12) $s_2=(2,0), t_1=(1,3), t_2=(3,2)$ (or $s_2=(2,1), t_1=(1,3), t_2=(3,2)$);  

(13) $s_2=(2,3), t_1=(1,3), t_2=(3,2)$ (or $s_2=(2,1), t_1=(1,3), t_2=(3,0)$);  

(14) $s_2=(2,3), t_1=(1,3), t_2=(3,2)$ (or $s_2=(2,3), t_1=(1,3), t_2=(3,1)$);  

(15) $s_2=(2,2), t_1=(1,3), t_2=(3,1)$.  \qed   

The following corollary is straightforward.

\begin{corollary}\label{base-main-e}{\bf.} Let $\{s_1,s_2\}$ and $\{t_1,t_2\}$ be any two sets of vertices in different partite sets of $BH_{2}$ and let $e$ be any edge of $BH_2$. Then there exist vertex-disjoint $s_1,t_1$-path and $s_2,t_2$-path of $BH_{2}-e$ that covers all vertices of it.
\end{corollary}

\noindent{\bf Remark.} Our aim is to guarantee that there exists a dimension $d\in\{0,1,2\}$ such that by dividing $BH_3$ into $B^i$ along dimension $d$ we can use Lemmas $\ref{base-main-e-f}$ and $\ref{base-main-e}$ as the induction basis of the main result. Let $F=\{f_0,f_1,f_2\}$ be a set of three edges of $BH_3$ and let $\{s_1,s_2\}$ and $\{t_1,t_2\}$ be any two sets of vertices in different partite sets of $BH_{3}$. If there exists a dimension $d\in\{0,1,2\}$ such that $|E_d\cap F|\geq2$, then $BH_3$ can be divided into $B^i$ ($0\leq i\leq 3$) along dimension $d$. Thus, $|E(B^i)\cap F|\leq1$ for each $i\in\{0,1,2,3\}$. So we assume that $E_j\cap F=\{f_i\}$ for each $j=0,1,2$. By Lemma \ref{common-nb}, $s_1$ and $s_2$ (or $t_1$ and $t_2$) have 0, 2 or $2n$ common neighbors.

If $s_1$ and $s_2$ (or $t_1$ and $t_2$) have no common neighbors, then we can safely divide $BH_3$ into $B^i$ ($0\leq i\leq 3$) along each dimension $d\in\{0,1,2\}$.

If $s_1$ and $s_2$ (or $t_1$ and $t_2$) have at least 2 common neighbors, we may assume that $x$ is one of the common neighbors of $s_1$ and $s_2$. If we divide $BH_3$ into $B^i$ ($0\leq i\leq 3$) along some dimension $d\in\{0,1,2\}$ such that $s_1$, $s_2$, $t_1$ and $t_2$ are in the same $B^i$, say $B^0$, and $F=F'$, where $F'$ is the set of edges incident to $x$ in $B^0$ (except $s_1x$ and $s_2x$). Furthermore, if $s_1$ and $s_2$ (or $t_1$ and $t_2$) have exact 2 common neighbors, then $s_1x$ and $s_2x$ are edges of different dimensions, then we can choose a dimension $d'\in\{0,1,2\}\setminus\{d\}$ such that by dividing $BH_3$ into $B^i$ ($0\leq i\leq 3$) along dimension $d'$, $s_1$ and $s_2$ (or $t_1$ and $t_2$) are not in the same $B^i$. If $s_1$ and $s_2$ have $6$ common neighbors, then $s_1x$ and $s_2x$ are edges of the same dimension, so we can divide $BH_3$ into $B^i$ ($0\leq i\leq 3$) along the dimension of the edges in $F'$.\qed


\begin{lemma}\label{exist-HP-abcd}{\bf.} Let $F$ be a set of edges of $BH_n$ ($n\geq3$) with $|F|=2n-3$. Given a dimension $k$ of $BH_n$ such that $|E_k\cap F|=\max\{|E_j\cap F||0\leq j\leq n-1\}$.
Let $B^i$, $0\leq i\leq 3$, be subgraphs of $BH_n$ obtained by splitting $BH_n$ along dimension $k$. Then there exists four vertices $a,c\in V_0$ and $b,d\in V_1$ of $B^i$ such that:%

(1) $a=p(c)$, $b=p(d)$, and $a,b,c$ and $d$ form a 4-cycle in $B^i$;

(2) there exists a $k$-dimension neighbor $a_{i+1}$ of $a$ and $c$ such that $e_k(a_{i+1})\cap F=\emptyset$;

(3) there exist two $k$-dimension neighbors $u_{i-1}$ and $v_{i-1}$ of $b$ and $d$ such that $e_k(b)\cap F=\emptyset, |e_k(d)\cap F|<2$ and $cd\not\in F$;

(4) there exists a neighbor $u$ of $b$ and $d$ in $B^i$ such that $|e_{j_1}(u)\cap F|<2$ for each $j_1\in\{0,1,\cdots,n-1\}$;

(5) there exists a longest path $P$ from $u$ to $a$ covering all vertices of $B^i-F$ but $b,c$ and $d$.

\end{lemma}

\noindent{\bf Proof.} We proceed the proof by induction on $n$. By the choice of $k$, we have $|E_k\cap F|=1$ or $|E_k\cap F|\geq2$ when $n=3$. It is easy to verify that conditions (1)-(5) hold after splitting $BH_3$ by dimension $k$. Thus, the induction basis holds. So we assume that the lemma is true for all integers $m$ with $3\leq m\leq n-1$. Next we consider $BH_n$.

\begin{figure}[h]
\centering
\includegraphics[width=60mm]{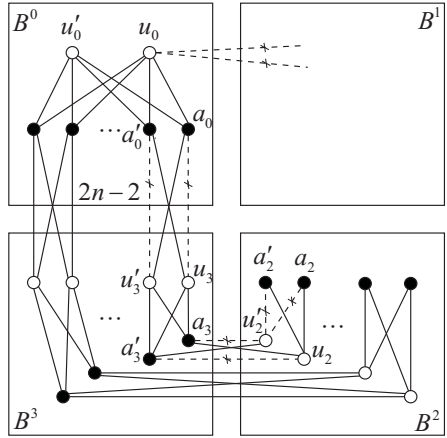}
\caption{Existence of $a,b,c$ and $d$ satisfying required conditions in Lemma \ref{exist-HP-abcd}.} \label{g_exist-abcd}
\end{figure}

Note that $|E_k\cap F|\geq2$ whenever $n\geq4$, suppose without loss of generality that $i=3$ and $k=n-1$. Since $|E_{n-1}\cap F|\geq2$, $|F\cap E(B^i)|\leq 2n-5$, $0\leq i\leq 3$. For each pair of vertices $u_0,u_0'\in V_0$ with $u_0=p(u_0')$ in $B^0$, there exist $2n-2$ common neighbors of them in $B^0$. Let $a_0$ and $a_0'$ be any two neighbors of $u_0$ and $u_0'$ with $a_0=p(a_0')$ in $B^0$. In addition, let $u_3$ and $u_3'$ be two ($n-1$)-dimension neighbors of $a_0$ and $a_0'$ and let $a_3,a_3'$ be two $k_1$-dimension neighbors of $u_3$ and $u_3'$ of $B^3$ for a given $k_1\in\{0,1,\cdots,n-2\}$. Accordingly, let $u_2$ and $u_2'$ be two ($n-1$)-dimension neighbors of $a_3$ and $a_3'$ and let $a_2$ and $a_2'$ be two $k_1$-dimension neighbors of $u_2$ and $u_2'$ in $B^2$. Thus, the subgraph induced by $\{a_0,a_0',u_3,u_3',a_3,a_3',u_2,u_2',a_2,a_2'\}$ is a twisted 4-cycle chain. If there exist at least two edges of $F$ in one of $\langle a_0',u_3,a_0,u_3',a_0'\rangle$, $\langle a_3',u_2,a_3,u_2',a_3'\rangle$ and $\langle a_2',u_2,a_2,u_2',a_2'\rangle$, then it may eliminate the choice of $a_3,a_3',u_3$ and $u_3'$ as $a,b,c$ and $d$ to satisfy conditions (1),(2) and (3) (see Fig. \ref{g_exist-abcd}). By arbitrary choice of $a_0$ and $a_0'$, if there exist no such $a,b,c$ and $d$ satisfying conditions (1),(2) and (3) for given $u_0$ and $u_0'$, we have $|F|=2\times(n-1)=2n-2>2n-3$, a contradiction.

On the other hand, $b$ and $d$ have $2n-2$ common neighbors (except $a$ and $c$) in $B^3$. Since $2\times(2n-2)>2\times(2n-4)>2n-3$ whenever $n\geq4$, there must exist a common neighbor $u$ of $b$ and $d$ satisfying condition (4). It remains to show that condition (5) holds.

By our assumption, $u,a,b,c,d\in V(B^3)$. Note that we have $|E(B^3)\cap F|\leq 2n-5$, our aim is to show that there exists a longest path $P$ from $u$ to $a$ covering all vertices of $B^3-F$ but $b,c$ and $d$. Let $k_2\in\{0,1,\cdots,n-2\}$ such that $|E_{k_2}\cap E(B^3)\cap F|\geq|E_{j}\cap E(B^3)\cap F|$ for each $j\in \{0,1,\cdots,n-2\}\setminus \{k_2\}$. We further divide each $B^i$ into $B_{n-2}^{i_1,i}$, $0\leq i_1\leq 3$, along dimension $k_2$. That is, $B_{n-2}^{i_1,i}\cong BH_{n-2}$ for each $i_1$ and $i$. Assume without loss of generality that $a,b,c,d\in V(B_{n-2}^{0,3})$. By Definition \ref{def1}, the graph induced by $V(B_{n-2}^{0,0})$, $V(B_{n-2}^{0,1})$, $V(B_{n-2}^{0,2})$ and $V(B_{n-2}^{0,3})$ is isomorphic to $BH_{n-1}$, for convenience, we denote it by $H$. Since $u$ is a neighbor of $b$ and $d$ in $B^3$, we assume without loss of generality that $u\in V(B_{n-2}^{0,3})$.

\begin{figure}[h]
\centering
\includegraphics[width=60mm]{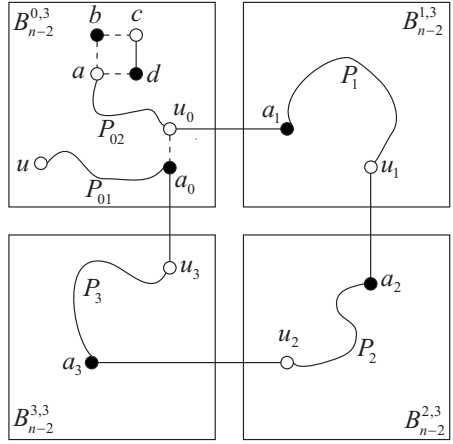}
\caption{Longest path from $u$ to $a$ covering all vertices of $B^3-F$ but $b,c$ and $d$.} \label{g_HP-abcd}
\end{figure}

By induction hypothesis, there exists a longest path $P_0$ from $u$ to $a$ covering all vertices of $B_{n-2}^{0,3}-F$ but $b,c$ and $d$. Since $l(P_0)=4^{n-2}-4$ and $(4^{n-2}-4)/2>2n-3$ whenever $n\geq4$ (any vertex $v$ on $P_0$ with $|e_{k_2}(v)\cap F|=2$ will eliminate the choice of two edges incident to $v$ on $P_0$), we can choose an edge $u_0a_0\in E(P_0)$ such that there exist two edges $u_0a_1,u_3a_0\not\in F$, where $a_1\in V(B_{n-2}^{1,3})$ and $u_3\in V(B_{n-2}^{3,3})$. Deleting $u_0a_0$ from $P_0$ will generate two vertex-disjoint paths $P_{01}$ and $P_{02}$, where $P_{01}$ connects $u$ to $a_0$ and $P_{02}$ connects $u_0$ to $a$. Let $u_1a_2$ and $u_2a_3$ be two fault-free $k_2$-dimension edges. By Lemma \ref{HP-fault}, there exist a fault-free Hamiltonian path $P_1$ of $B_{n-2}^{1,3}$ from $u_1$ to $a_1$, a fault-free Hamiltonian path $P_2$ of $B_{n-2}^{2,3}$ from $u_2$ to $a_2$, and a fault-free Hamiltonian path $P_3$ of $B_{n-2}^{3,3}$ from $u_3$ to $a_3$. Hence, $\langle u,P_{01},a_0,u_3,P_3,a_3,u_2,P_2,a_2,u_1,P_1,a_1,u_0,P_{02},a\rangle$ is the path required (see Fig. \ref{g_HP-abcd}).

This completes the proof. \qed

\begin{lemma}\label{exist-t_1_abc-bh2-f2}{\bf.} Let $F=\{e,f\}$ be any two edges of $BH_2$ with $e\in E_0$ and $f\in E_1$. In addition, let $t_1,t_2\in V_1$ be two arbitrary vertices. Then there exist two pairs of vertices in $V_0$ differing only from inner index respectively, suppose without loss of generality that $a$ and $c$ is such a pair with $a=p(c)$, such that: (1) there exists a vertex $u\in V_0$ of $BH_2$ with $u\neq a,c$; (2) there exist two vertex-disjoint paths $P$ and $Q$ of $BH_2-F$ cover all vertices of it, where $P$ connects $u$ to $t_2$, and $Q$ connects $c$ to $t_1$ and $\langle c,b,a\rangle$ is a subpath of $Q$.
\end{lemma}
\noindent{\bf Proof.} By vertex-transitivity of $BH_2$, we may assume that $t_1=(1,0)$. Since $e\in E_0$ and $f\in E_1$, $e$ and $f$ lie in different twisted 4-cycles of $BH_2$. Our aim is to find two pairs of vertices differing only from inner index respectively and satisfying conditions (1) and (2). There are three essentially different positions of $t_2$.

\noindent{\bf Case 1.} $t_2=(3,0)$. We further deal with the following cases.

\noindent{\bf Case 1.1.} $e$ and $f$ lie in consecutive twisted 4-cycles.

\noindent{\bf Case 1.1.1.} $e$ and $f$ are nonadjacent. We may assume that $e=(0,0)(1,0)$ and $f=(2,0)(3,1)$. If $a=(2,0)$, $c=(0,0)$ and $u=(2,1)$, then $P^{-1}=\langle T_3((3,0);(0,1)),(3,1),(2,1)\rangle$ and $Q=\langle (0,0),(1,1),(2,0),(1,0)\rangle$ are the paths required.

If $a=(0,1)$, $c=(2,1)$ and $u=(2,2)$, then $P=\langle$(2,2),(3,2),(0,2),(1,3),(0,3), (3,3),(2,3),(3,0)$\rangle$ and $Q=\langle (2,1),(1,2),(0,1),(3,1)$, $(0,0),(1,1),(2,0),(1,0)\rangle$ are the paths required.

\noindent{\bf Case 1.1.2.} $e$ and $f$ are adjacent. There are two essential relative positions of $e$ and $f$, we further deal with the following cases.

\noindent{\bf Case 1.1.2.1.} $e=(0,0)(1,0)$ and $f=(0,0)(1,1)$. The proof is similar to that of Case 1.1.1, we omit it.

\noindent{\bf Case 1.1.2.2.} $e=(0,0)(1,1)$ and $f=(1,1)(0,1)$. If $a=(2,1)$, $c=(0,1)$ and $u=(2,2)$, then $P=\langle$(2,2),(3,2),(0,2),(1,3),(0,3),(3,3),(2,3),(3,0)$\rangle$ and $Q=\langle (0,1),(1,2),(2,1),(1,1)$, $(2,0),(3,1),(0,0),(1,0)\rangle$ are the paths required.

If $a=(0,2)$, $c=(2,2)$ and $u=(2,3)$, then $P=\langle (2,3),(1,3),(0,3),(1,0)\rangle$ and $Q=\langle (2,2),(3,3),(0,2),(1,2)$, $(0,1),(3,2),(2,1),(1,1),(2,0),(3,1),(0,0),(1,0)\rangle$ are the paths required.

\noindent{\bf Case 1.2.} $e$ and $f$ lie in inconsecutive twisted 4-cycles. Obviously, $BH_2$ can be decomposed into four edge-disjoint 8-cycles according to ring-like layout. By Lemma \ref{T_uvxy}, each pair of vertices in $V_0$ differing only from the inner index can be chosen as $a$ and $c$ such that there exist two vertex-disjoint paths $P$ and $Q$ of $BH_2-F$ cover all vertices of it, where $P$ connects $u$ to $t_2$, and $Q$ connects $c$ to $t_1$ and $\langle c,b,a\rangle$ is a subpath of $Q$.

\noindent{\bf Case 2.} $t_2=(3,3)$. We further deal with the following cases.

\noindent{\bf Case 2.1.} $|F\cap T_0(t_1,(3,0);(1,3),t_2)|=2$. By Lemma \ref{T_uvxy}, there exist two edge-disjoint 2-paths $P_1$ and $Q_1$ cover all vertices of $T_0(t_1,(3,0);(1,3),t_2)$, where $P_1$ connects (3,0) to $t_2$ and $Q_1$ connects (1,3) to $t_1$. There are two pairs of vertices can be chosen as $a$ and $c$: (1) $a=(0,2)$ and $c=(2,2)$; (2) $a=(0,1)$ and $c=(2,1)$.

If $a=(0,2)$ and $c=(2,2)$, let $u=(2,1)$, then $P=\langle (2,1),(1,2),(0,1),(3,1)$, $(0,0),(1,1),(2,0),(3,0),P_1,(3,3)\rangle$ and $Q=\langle (2,2),(3,2),(0,2),(1,3),Q_1,(1,0)\rangle$ are the paths required.

If $a=(0,1)$ and $c=(2,1)$, let $u=(2,0)$, then $P=\langle (2,0),(3,1)$, $(0,0),(3,0),P_1$, (3,3)$\rangle$ and $Q=\langle (2,1),(1,1),(0,1),(1,2),(0,2),(3,2),(2,2),(1,3)$, $Q_1,(1,0)\rangle$ are the paths required.

\noindent{\bf Case 2.2.} $|F\cap T_0(t_1,(3,0);(1,3),t_2)|=1$ or $|F\cap T_0(t_1,(3,0);(1,3),t_2)|=0$. The proof is similar to that of Case 2.1, we omit it.

\noindent{\bf Case 3.} $t_2=(3,2)$. The proof is similar to that of Case 2, we omit it. \qed

By above lemma, it is not hard to obtain the following corollary.

\begin{corollary}\label{exist-t_1_abc-bh2-f1}{\bf.} Let $e$ be an edge of $BH_2$. In addition, let $t_1,t_2\in V_1$ be two arbitrary vertices of $BH_2$. Then there exist at least two pairs of vertices in $V_0$ differing only from inner index respectively, suppose without loss of generality that $a$ and $c$ is such a pair with $a=p(c)$, such that: (1) there exists a vertex $u\in V_0$ of $BH_2$ with $u\neq a,c$; (2) there exist two vertex-disjoint paths $P$ and $Q$ of $BH_2-e$ cover all vertices of it, where $P$ connects $u$ to $t_2$, and $Q$ connects $c$ to $t_1$ and $\langle c,b,a\rangle$ is a subpath of $Q$.
\end{corollary}

\begin{lemma}\label{exist-HP-t_1_abc}{\bf.} Let $F$ be a set of edges of $BH_n$ with $|F|=2n-3$ ($n\geq3$). Given a dimension $k$ of $BH_n$ such that $|E_k\cap F|\geq |E_j\cap F|$ for each $j\in\{0,1,\cdots,n-1\}\setminus\{k\}$.
Let $B^i$, $0\leq i\leq 3$, be subgraphs of $BH_n$ obtained by splitting $BH_n$ along dimension $k$. In addition, let $t_1,t_2\in V_1$ be two arbitrary vertices in $B^i$ such that $t_1\neq t_2$. Then, there exist four vertices $u,a,c\in V_0$ and $b\in V_1$ of $B^i$ with $a=p(c)$ such that:

(1) there exists a $k$-dimension neighbor $a_{i+1}$ of $a$ and $c$ such that $e_k(a_{i+1})\cap F=\emptyset$ and there exists a $k$-dimension neighbor $u_{i-1}$ of $b$ such that $|e_k(b)\cap F|<2$, where $b$ ($b\neq t_1,t_2$) is a common neighbor of $a$ and $c$;

(2) for each $j_1\in\{0,1,\cdots,n-1\}$, $|e_{j_1}(u)\cap F|<2$;

(3) there exist two vertex-disjoint paths $P$ and $Q$ of $B^i-F$ cover all vertices of it, where $P$ connects $u$ to $t_2$, and $Q$ connects $c$ to $t_1$ and $\langle c,b,a\rangle$ is a subpath of $Q$.
\end{lemma}

\noindent{\bf Proof.} We proceed the proof by induction on $n$. Firstly, we shall show that the lemma is true when $n=3$. Suppose without loss of generality that $i=3$ and $k=2$, that is, $t_1,t_2\in V(B^3)$. Since $|E_{2}\cap F|\geq1$, $|F\cap E(B^i)|\leq 2$ for $0\leq i\leq 3$. If $|E_{2}\cap F|=1$, then it follows from Lemma \ref{exist-t_1_abc-bh2-f2} that the lemma is true. If $|E_{2}\cap F|\geq2$, it follow from Lemma \ref{exist-t_1_abc-bh2-f1} that the lemma is also true. Thus, the induction basis holds. So we assume that the lemma is true for all integers $m$ with $3\leq m\leq n-1$. Next we consider $BH_n$.

Obviously, we have $|E_{k}\cap F|\geq2$ whenever $n\geq4$. We may assume that $i=3$ and $k=n-1$. So we obtain four subgraphs $B^i$, $0\leq i\leq 3$, by splitting $BH_n$ along dimension $n-1$. Accordingly, by our assumption, $t_1,t_2\in V(B^3)$. Thus, we have $|E(B^3)\cap F|\leq 2n-5$. Our aim is to show that there exist four vertices $u,a,c\in V_0$ and $b\in V_1$ of $B^3$ with $a=p(c)$ satisfying conditions (1)-(3). Let $k_1\in\{0,1,\cdots,n-2\}$ such that $|E_{k_1}\cap E(B^3)\cap F|\geq|E_{j}\cap E(B^3)\cap F|$ for each $j\in \{0,1,\cdots,n-2\}\setminus \{k_1\}$. We further divide each $B^i$ into $B_{n-2}^{i_1,i}$, $0\leq i_1\leq 3$, along dimension $k_1$. That is, $B_{n-2}^{i_1,i}\cong BH_{n-2}$ for each $i_1$ and $i$. Assume without loss of generality that $t_1\in V(B_{n-2}^{0,3})$. By Definition \ref{def1}, the graph induced by $V(B_{n-2}^{0,0})$, $V(B_{n-2}^{0,1})$, $V(B_{n-2}^{0,2})$ and $V(B_{n-2}^{0,3})$ is isomorphic to $BH_{n-1}$, for convenience, we denote it by $H$. There are four relative positions of $t_2$ in $B^3$, so we consider the following conditions.

If $t_2\in V(B_{n-2}^{0,3})$. By induction hypothesis, there exist four vertices $u,a,c\in V_0$ and $b\in V_1$ of $B_{n-2}^{0,3}$ with $a=p(c)$ satisfying conditions (1),(2) in $H$. Moreover, there exist two vertex-disjoint paths $P_0$ and $Q$ of $B_{n-2}^{0,3}-F$ cover all vertices of it, where $P_0$ connects $u$ to $t_2$, and $Q$ connects $c$ to $t_1$ and $\langle c,b,a\rangle$ is a subpath of $Q$. Since $l(P_0)+l(Q)=4^{n-2}-2$, it is obvious that there exists an edge on $P_0$ or $Q$, say $u_0a_0\in E(P_0)$, such that $u_0a_1,u_3a_0\not\in F$, where $u_0a_1$ and $u_3a_0$ are $k_1$-dimension edges. Thus, deleting $u_0a_0$ from $P_0$ will generate two vertex-disjoint paths $P_{01}$ and $P_{02}$, where $P_{01}$ connects $u$ to $a_0$ and $P_{02}$ connects $u_0$ to $t_2$. By Lemma \ref{crossing-edges}, there must exist two $k_1$-dimension fault-free edges $u_1a_2$ and $u_2a_3$, where $u_1\in V(B_{n-2}^{1,3})$, $u_2,a_2\in V(B_{n-2}^{2,3})$ and $a_3\in V(B_{n-2}^{3,3})$. By Lemma \ref{HP-fault}, there exist a fault-free Hamiltonian path $P_1$ of $B_{n-2}^{1,3}-F$ from $u_1$ to $a_1$, a fault-free Hamiltonian path $P_2$ of $B_{n-2}^{2,3}-F$ from $u_2$ to $a_2$, and a fault-free Hamiltonian path $P_3$ of $B_{n-2}^{3,3}-F$ from $u_3$ to $a_3$. Hence, $P=\langle u,P_{01},a_0,u_3,P_{3},a_3,u_2,P_{2},a_2,u_1,P_{1},a_1,u_0,P_{02},t_2\rangle$ and $Q$ are paths satisfying condition (3) in $BH_n$.

If $t_2\in V(B_{n-2}^{1,3})$. Obviously, there exists a vertex $u\in V(B_{n-2}^{1,3})$ such that $|e_{j_1}(u)\cap F|<2$ for each $j_1\in\{0,1,\cdots,n-1\}$. By Lemma \ref{HP-fault}, there exists a fault-free Hamiltonian path $P_1$ of $B_{n-2}^{1,3}-F$ from $u$ to $t_2$. Since $l(P_1)=4^{n-2}-1$, there must exist an edge $u_1a_1\in E(P_1)$ such that $|e_{k_1}(u_1)\cap F|<2$ and $|e_{k_1}(a_1)\cap F|<2$. So let $u_1a_2$ and $u'a_1$ be two fault-free $k_1$-dimension edges. Additionally, deleting $u_1a_1$ from $P_1$ will generate two vertex-disjoint paths $P_{11}$ and $P_{12}$, where $P_{11}$ connects $u$ to $a_1$ and $P_{12}$ connects $u_1$ to $t_2$. By induction hypothesis, there exists four vertices $a,c\in V_0$ and $a_0,b\in V_1$ of $B_{n-2}^{0,3}$ with $a=p(c)$ such that: $a,b$ and $c$ satisfy condition (1) and $a_0$ satisfies condition (2) in $H$. Moreover, there exist two vertex-disjoint paths $P_0$ and $Q$ of $B_{n-2}^{0,3}-F$ cover all vertices of it, where $P_0$ connects $u'$ to $a_0$, and $Q$ connects $c$ to $t_1$ and $\langle c,b,a\rangle$ is a subpath of $Q$. Obviously, there exist two $k_1$-dimension fault-free edges $u_2a_3$ and $u_3a_0$, where $u_2\in V(B_{n-2}^{2,3})$ and $u_3,a_3\in V(B_{n-2}^{3,3})$. By Lemma \ref{HP-fault}, there exist a fault free Hamiltonian path $P_2$ of $B_{n-2}^{2,3}-F$ from $u_2$ to $a_2$, and a fault free Hamiltonian path $P_3$ of $B_{n-2}^{3,3}-F$ from $u_3$ to $a_3$. Hence, $P=\langle u,P_{11},a_1,u',P_0,a_0,u_3,P_{3},a_3,u_2,P_{2},a_2,u_1,P_{12},t_2\rangle$ and $Q$ are paths satisfying condition (3) in $BH_n$.

If $t_2\in V(B_{n-2}^{2,3})$. Obviously, there exists a vertex $u\in V_0$ in $B_{n-2}^{2,3}$ such that $|e_{j_1}(u)\cap F|<2$ for each $j_1\in\{0,1,\cdots,n-1\}$. By Lemma \ref{HP-fault}, there exists a fault-free Hamiltonian path $P_2$ of $B_{n-2}^{2,3}-F$ from $u$ to $t_2$. Similarly, there must exist an edge $u_2a_2\in E(P_2)$ such that $|e_{k_1}(u_2)\cap F|<2$ and $|e_{k_1}(a_2)\cap F|<2$. So let $u_1a_2$ and $u_2a_3$ be two fault-free $k_1$-dimension edges. Additionally, deleting $u_2a_2$ from $P_2$ will generate two vertex-disjoint paths $P_{21}$ and $P_{22}$, where $P_{21}$ connects $u$ to $a_2$ and $P_{22}$ connects $u_2$ to $t_2$. Let $a_0\in V(B_{n-2}^{0,3})$ be a vertex such that $|e_{k_1}(a_0)\cap F|<2$. By induction hypothesis, there exists four vertices $u_0,a,c\in V_0$ and $b\in V_1$ of $B_{n-2}^{0,3}$ with $a=p(c)$ such that: $a,b$ and $c$ satisfy condition (1) and $u_0$ satisfies condition (2) in $H$. Moreover, there exist two vertex-disjoint paths $P_0$ and $Q$ of $B_{n-2}^{0,3}-F$ cover all vertices of it, where $P_0$ connects $u_0$ to $a_0$, and $Q$ connects $c$ to $t_1$ and $\langle c,b,a\rangle$ is a subpath of $Q$. Obviously, there exist two fault-free $k_1$-dimension edges $u_0a_1$ and $u_3a_0$, where $u_3\in V(B_{n-2}^{3,3})$ and $a_1\in V(B_{n-2}^{1,3})$. By Lemma \ref{HP-fault}, there exist a fault-free Hamiltonian path $P_1$ of $B_{n-2}^{1,3}-F$ from $u_1$ to $a_1$, and a fault free Hamiltonian path $P_3$ of $B_{n-2}^{3,3}-F$ from $u_3$ to $a_3$. Hence, $P=\langle u,P_{21},a_2,u_1,P_1,a_1,u_0,P_{0},a_0,u_3,P_{3},a_3,u_2,P_{22},t_2\rangle$ and $Q$ are paths satisfying condition (3) in $BH_n$.

If $t_2\in V(B_{n-2}^{3,3})$. Similarly, there exists a vertex $u\in V(B_{n-2}^{2,3})$ such that $|e_{j_1}(u)\cap F|<2$ for each $j_1\in\{0,1,\cdots,n-1\}$. Let $a_0\in V(B_{n-2}^{0,3})$ be a vertex such that $|e_{k_1}(a_0)\cap F|<2$. By induction hypothesis, there exists four vertices $u_0,a,c\in V_0$ and $b\in V_1$ of $B_{n-2}^{0,3}$ with $a=p(c)$ such that: $a,b$ and $c$ satisfy condition (1) and $u_0$ satisfies condition (2) in $H$. Moreover, there exist two vertex-disjoint paths $P_0$ and $Q$ of $B_{n-2}^{0,3}-F$ cover all vertices of it, where $P_0$ connects $u_0$ to $a_0$, and $Q$ connects $c$ to $t_1$ and $\langle c,b,a\rangle$ is a subpath of $Q$. So there exist three fault-free $k_1$-dimension edges $u_0a_1$, $u_1a_2$ and $u_3a_0$, where $u_1,a_1\in V(B_{n-2}^{1,3})$, $a_2\in V(B_{n-2}^{2,3})$ and $u_3\in V(B_{n-2}^{3,3})$. By Lemma \ref{HP-fault}, there exist a fault-free Hamiltonian path $P_1$ of $B_{n-2}^{1,3}-F$ from $u_1$ to $a_1$, a fault-free Hamiltonian path $P_2$ of $B_{n-2}^{2,3}-F$ from $u$ to $a_2$ and a fault-free Hamiltonian path $P_3$ of $B_{n-2}^{3,3}-F$ from $u_3$ to $t_2$. Hence, $P=\langle u,P_{2},a_2,u_1,P_1,a_1,u_0,P_{0},a_0,u_3,P_{3},t_2\rangle$ and $Q$ are paths satisfying condition (3) in $BH_n$.

This completes the proof. \qed

Now we are ready to state the main result of this paper.

\begin{theorem}{\bf.} Let $F$ be a set of edges with $|F|\leq 2n-3$ and let $\{s_1,s_2\}$ and $\{t_1,t_2\}$ be two sets of vertices in different partite sets of $BH_{n}$ for $n\geq2$. Then $BH_{n}-F$ contains vertex-disjoint $s_1,t_1$-path and $s_2,t_2$-path that covers all vertices of it.
\end{theorem}

\noindent{\bf Proof.} We proceed the proof by induction on $n$. By Lemma \ref{base-main-e}, the statement is true for $BH_2$. For $n=3$, we have characterized how to divide $BH_3$ by some dimension $k\in\{0,1,2\}$ in Remark. Assume that the statement holds for $BH_{n-1}$ with $n\geq3$. Next we consider $BH_n$. Since $|F|\leq 2n-3$, by Pigeonhole Principle, there must exist a dimension $d\in\{0,1,\cdots,n-1\}$ such that $|E_d\cap F|\geq2$ whenever $n\geq4$. Thus, $|E(B^i)\cap F|\leq 2n-5$, $i\in\{0,1,2,3\}$ (we can also use Lemma \ref{base-main-e-f} as induction basis when $n=3$). Suppose without loss of generality that $d=n-1$. So we divide $BH_n$ into four subcubes $B^{i}$ ($i\in\{0,1,2,3\}$) by deleting $E_{n-1}$. By Lemma \ref{transitive}, $BH_n$ is vertex transitive, we may assume that $s_1\in V(B^0)$ and $|V(B^0)\cap \{s_2,t_1,t_2\}|\geq |V(B^j)\cap \{s_2,t_1,t_2\}|$ for $j\in\{1,2,3\}$. We consider the following cases.

\vskip 0.05 in

\noindent{\bf Case 1.} $|V(B^0)\cap \{s_2,t_1,t_2\}|=0$. We further deal with the following cases.

\begin{figure}
\begin{minipage}[t]{0.5\linewidth}
\centering
\includegraphics[width=60mm]{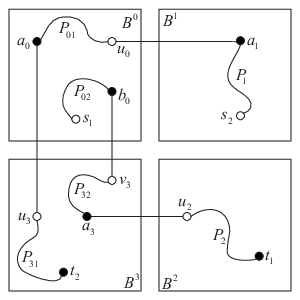}
\caption{Illustration of Cases 1.1.} \label{g3}
\end{minipage}
\begin{minipage}[t]{0.5\linewidth}
\centering
\includegraphics[width=60mm]{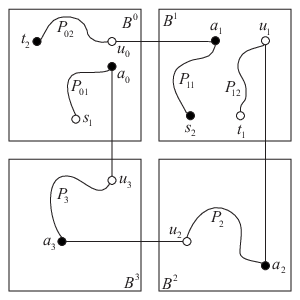}
\caption{Illustration of Case 2.1.2.} \label{g4}
\end{minipage}
\end{figure}

\noindent{\bf Case 1.1.} $s_2\in V(B^1)$, $t_1\in V(B^2)$ and $t_2\in V(B^3)$. Since $4^{n-1}\geq 2n-3$ whenever $n\geq3$, there always exists a fault-free edge $u_3a_0\in E_{3,0}$. In addition, there exists a fault-free edge $v_3b_0\in E_{3,0}$ such that $u_3\neq v_3$ and $b_0\neq a_0$ (Let $u$ and $v$ be two vertices with distance two and let $w$ be a common neighbor of them in $B^i$. We denote the set of edges incident to $w$, except $uw$ and $vw$, in $B^i$ by $A$, then $|A|=2n-4$. Note that $|F^i|\leq2n-4$ when $n=3$ and $|F^i|\leq2n-5$ when $n>3$, there may exist at most one pair of vertices $u$ and $v$ such that $A\subseteq F$ only if $n=3$. It is easy to choose $u$ and $v$ to avoid this situation, so we do not mention this condition in the following proof.). Similarly, there exist two fault-free edges $u_0a_1\in E_{0,1}$ and $u_2a_3\in E_{2,3}$ such that $u_0\neq s_1$ and $a_3\neq t_2$. By Lemma \ref{HP-fault}, there exist a fault-free Hamiltonian path $P_1$ of $B^1-F$ from $s_2$ to $a_1$, and a fault-free Hamiltonian path $P_2$ of $B^2-F$ from $u_2$ to $t_1$. By induction hypothesis, there exist two vertex-disjoint paths $P_{01}$ and $P_{02}$ cover all vertices of $B^0-F$, where $P_{01}$ connects $u_0$ to $a_0$ and $P_{02}$ connects $s_1$ to $b_0$; there exist two vertex-disjoint paths $P_{31}$ and $P_{32}$ cover all vertices of $B^3-F$, where $P_{31}$ connects $u_3$ to $t_2$ and $P_{32}$ connects $v_3$ to $a_3$. Hence, $\langle s_1,P_{02},b_0,v_3,P_{32},a_3,u_2,P_2,t_1\rangle$ and $\langle s_2,P_1,a_1, u_0,P_{01},a_0,u_3,P_{31},t_2\rangle$ are two vertex-disjoint paths required (see Fig. \ref{g3}). 

\noindent{\bf Case 1.2.} $s_2\in V(B^1)$, $t_1\in V(B^3)$ and $t_2\in V(B^2)$. There always exist two edges $u_3a_0,v_3b_0\in E_{3,0}$ such that $u_3\neq v_3$ and $a_0\neq b_0$. Similarly, there exist an edge $u_0a_1\in E_{0,1}$ such that $u_0\neq s_1$, and an edge $u_2a_3\in E_{2,3}$ such that $a_3\neq t_1$. By Lemma \ref{HP-fault}, there exist a fault-free Hamiltonian path $P_1$ of $B^1-F$ from $s_2$ to $a_1$, and a fault-free Hamiltonian path $P_2$ of $B^2-F$ from $u_2$ to $t_2$. By induction hypothesis, there exist two vertex-disjoint paths $P_{31}$ and $P_{32}$ cover all vertices of $B^3-F$, where $P_{31}$ connects $v_3$ to $t_1$ and $P_{32}$ connects $u_3$ to $a_3$; there exist two vertex-disjoint paths $P_{01}$ and $P_{02}$ cover all vertices of $B^0-F$, where $P_{01}$ connects $u_0$ to $a_0$ and $P_{02}$ connects $s_1$ to $b_0$. Hence, $\langle s_1,P_{02},b_0,v_3,P_{31},t_1\rangle$ and $\langle s_2,P_1,a_1, u_0,P_{01},a_0,u_3,P_{32},a_3,u_2,P_{2},t_2\rangle$ are two vertex-disjoint paths required.

\noindent{\bf Case 1.3.} $s_2\in V(B^2)$, $t_1\in V(B^1)$ and $t_2\in V(B^3)$.
There always exist two edges $u_3a_0,v_3b_0\in E_{3,0}$ such that $u_3\neq v_3$ and $a_0\neq b_0$, and two edges $u_1a_2,v_1b_2\in E_{1,2}$ such that $u_1\neq v_1$ and $a_2\neq b_2$.
Similarly, there exist an edge $u_0a_1\in E_{0,1}$ such that $u_0\neq s_1$ and $a_1\neq t_1$, and an edge $u_2a_3\in E_{2,3}$ such that $u_2\neq s_2$ and $a_3\neq t_2$. By induction hypothesis, there exist two vertex-disjoint paths $P_{01}$ and $P_{02}$ cover all vertices of $B^0-F$, where $P_{01}$ connects $u_0$ to $a_0$ and $P_{02}$ connects $s_1$ to $b_0$; there exist two vertex-disjoint paths $P_{11}$ and $P_{12}$ cover all vertices of $B^1-F$, where $P_{11}$ connects $v_1$ to $t_1$ and $P_{12}$ connects $u_1$ to $a_1$; there exist two vertex-disjoint paths $P_{21}$ and $P_{22}$ cover all vertices of $B^2-F$, where $P_{21}$ connects $u_2$ to $b_2$ and $P_{22}$ connects $s_2$ to $a_2$; there exist two vertex-disjoint paths $P_{31}$ and $P_{32}$ cover all vertices of $B^3-F$, where $P_{31}$ connects $v_3$ to $a_3$ and $P_{32}$ connects $u_3$ to $t_2$. Hence, $\langle s_1,P_{02},b_0,v_3,P_{31},a_3,u_2,P_{21},b_2,v_1,P_{11},t_1\rangle$ and $\langle s_2,P_{22},a_2,u_1,P_{12},a_1,u_0,P_{01},a_0,u_3,P_{32},t_2\rangle$ are two vertex-disjoint paths required.

\noindent{\bf Case 1.4.} $s_2\in V(B^2)$, $t_1\in V(B^3)$ and $t_2\in V(B^1)$. Obviously, there exist two non-faulty edges $u_3a_0\in E_{3,0}$ and $u_1a_2\in E_{1,2}$. By Lemma \ref{HP-fault}, there exist a fault-free Hamiltonian path $P_0$ of $B^0-F$ from $s_1$ to $a_0$, a fault-free Hamiltonian path $P_1$ of $B^1-F$ from $u_1$ to $t_2$, a fault-free Hamiltonian path $P_2$ of $B^2-F$ from $s_2$ to $a_2$, and a fault-free Hamiltonian path $P_3$ of $B^3-F$ from $u_3$ to $t_1$. Hence, $\langle s_1,P_{0},a_0,u_3,P_{3},t_1\rangle$ and $\langle s_2,P_{2},a_2,u_1,P_{1},t_2\rangle$ are two vertex-disjoint paths required.


\noindent{\bf Case 2.} $|V(B^0)\cap \{s_2,t_1,t_2\}|=1$. We further deal with the following cases.

\noindent{\bf Case 2.1.} For some $j\in \{0,1,2\}$, $|V(B^j)\cap \{s_2,t_1,t_2\}|=2$. %

\noindent{\bf Case 2.1.1.} $t_1\in V(B^0)$ and $s_2,t_2\in V(B^1)$. By Lemma \ref{HP-fault}, there exists a fault-free Hamiltonian path $P_0$ of $B^0-F$ from $s_1$ to $t_1$. Since $4^{n-1}-3\geq 2(2n-3)$ whenever $n\geq3$ and any vertex incident to two faulty $(n-1)$-dimension edges will eliminate the choice of two edges on $P_0$, we can choose an edge $u_0a_0\in E(P_0)$ such that there exist two non-faulty edges $u_0a_1\in E_{0,1}$ and $u_3a_0\in E_{3,0}$. Deleting $u_0a_0$ from $P_0$ will give rise to two disjoint paths $P_{01}$ and $P_{02}$, where $P_{01}$ connects $s_1$ to $u_0$ and $P_{02}$ connects $a_0$ to $t_1$. Additionally, there exist a fault-free edge $u_1a_2\in E_{1,2}$ such that $u_1\neq s_2$, and an edge $u_2a_3\in E_{2,3}$. By induction hypothesis, there exist two vertex-disjoint paths $P_{11}$ and $P_{12}$ cover all vertices of $B^1-F$, where $P_{11}$ connects $a_1$ to $u_1$ and $P_{12}$ connects $s_2$ to $t_2$. Moreover, there exist a fault-free Hamiltonian path $P_2$ of $B^2-F$ from $a_2$ to $u_2$, and a fault-free Hamiltonian path $P_3$ of $B^3-F$ from $a_3$ to $u_3$. Hence, $\langle s_1,P_{01},u_0,a_1,P_{11},u_1,a_2,P_{2},u_2,a_3,P_3,u_3,a_0,P_{02},t_1\rangle$ and $\langle s_2,P_{12},t_2\rangle$ are two vertex-disjoint paths required.

\noindent{\bf Case 2.1.2.} $t_2\in V(B^0)$ and $s_2,t_1\in V(B^1)$. There exist fault-free edges $u_0a_1\in E_{0,1}$ such that $u_0\neq s_1$ and $a_1\neq t_1$, $u_1a_2\in E_{1,2}$ such that $u_1\neq s_2$, $u_2a_3\in E_{2,3}$ and $u_3a_0\in E_{3,0}$ such that $a_0\neq t_2$. By induction hypothesis, there exist two vertex-disjoint paths $P_{01}$ and $P_{02}$ cover all vertices of $B^0-F$, where $P_{01}$ connects $s_1$ to $a_0$ and $P_{02}$ connects $u_0$ to $t_2$; there exist two vertex-disjoint paths $P_{11}$ and $P_{12}$ cover all vertices of $B^1-F$, where $P_{11}$ connects $u_1$ to $t_1$ and $P_{12}$ connects $s_2$ to $a_1$. By Lemma \ref{HP-fault}, there exist a fault-free Hamiltonian path $P_2$ of $B^2-F$ from $u_2$ to $a_2$, and a fault-free Hamiltonian path $P_3$ of $B^3-F$ from $u_3$ to $a_3$. Hence, $\langle s_1,P_{01},a_0,u_3,P_{3},a_3,u_2,P_{2},a_2,u_1,P_{11},t_1\rangle$ and $\langle s_2,P_{12},a_1,u_0,P_{02},t_2\rangle$ are two vertex-disjoint paths required (see Fig. \ref{g4}).

\noindent{\bf Case 2.1.3.} $t_1\in V(B^0)$ and $s_2,t_2\in V(B^2)$. By Lemma \ref{HP-fault}, there exists a fault-free Hamiltonian path $P_0$ of $B^0-F$ from $s_1$ to $t_1$. We can choose an edge $u_0a_0\in E(P_0)$ such that there exist two fault-free edges $u_0a_1\in E_{0,1}$ and $u_3a_0\in E_{3,0}$. Deleting $u_0a_0$ from $P_0$ will give rise to two disjoint paths $P_{01}$ and $P_{02}$, where $P_{01}$ connects $s_1$ to $u_0$ and $P_{02}$ connects $a_0$ to $t_1$. There exist a fault-free edge $u_1a_2\in E_{1,2}$ such that $a_2\neq t_2$, and a fault-free edge $u_2a_3\in E_{2,3}$ such that $u_2\neq s_2$. By induction hypothesis, there exist two vertex-disjoint paths $P_{21}$ and $P_{22}$ cover all vertices of $B^2-F$, where $P_{21}$ connects $a_2$ to $u_2$ and $P_{22}$ connects $s_2$ to $t_2$. By Lemma \ref{HP-fault}, there exist a fault-free Hamiltonian path $P_1$ of $B^1-F$ from $a_1$ to $u_1$, and a fault-free Hamiltonian path $P_3$ of $B^3-F$ from $a_3$ to $u_3$. Hence, $\langle s_1,P_{01},u_0,a_1,P_{1},u_1,a_2,P_{21},u_2,a_3,P_{3},u_3,a_0,P_{02},t_1\rangle$ and $\langle s_2,P_{22},t_2\rangle$ are two vertex-disjoint paths required.

\noindent{\bf Case 2.1.4.} $t_2\in V(B^0)$ and $s_2,t_1\in V(B^2)$. There exist fault-free edges $u_0a_1\in E_{0,1}$ such that $u_0\neq s_1$, $u_1a_2\in E_{1,2}$ such that $a_2\neq t_1$, $u_2a_3\in E_{2,3}$ such that $u_2\neq s_2$ and $u_3a_0\in E_{3,0}$ such that $a_0\neq t_2$. By induction hypothesis, there exist two vertex-disjoint paths $P_{01}$ and $P_{02}$ cover all vertices of $B^0-F$, where $P_{01}$ connects $s_1$ to $a_0$ and $P_{02}$ connects $u_0$ to $t_2$; there exist two vertex-disjoint paths $P_{21}$ and $P_{22}$ cover all vertices of $B^2-F$, where $P_{21}$ connects $u_2$ to $t_1$ and $P_{22}$ connects $s_2$ to $a_2$. By Lemma \ref{HP-fault}, there exist a fault-free Hamiltonian path $P_1$ of $B^1-F$ from $u_1$ to $a_1$, and a fault-free Hamiltonian path $P_3$ of $B^3-F$ from $u_3$ to $a_3$. Hence, $\langle s_1,P_{01},a_0,u_3,P_{3},a_3,u_2,P_{21},t_1\rangle$ and $\langle s_2,P_{22},a_2,u_1,P_{1},a_1,u_0,P_{02},t_2\rangle$ are two vertex-disjoint paths required.

\noindent{\bf Case 2.1.5.} $s_2\in V(B^0)$ and $t_1,t_2\in V(B^1)$. There always exist two fault-free edges $u_3a_0,v_3b_0\in E_{3,0}$ such that $u_3\neq v_3$ and $a_0\neq b_0$, two fault-free edges $u_1a_2,v_1b_2\in E_{1,2}$ such that $u_1\neq v_1$ and $a_2\neq b_2$, and two fault-free edges $u_2a_3,v_2b_3\in E_{2,3}$ such that $u_2\neq v_2$ and $a_3\neq b_3$. By induction hypothesis, there exist two vertex-disjoint paths $P_{01}$ and $P_{02}$ cover all vertices of $B^0-F$, where $P_{01}$ connects $s_1$ to $a_0$ and $P_{02}$ connects $s_2$ to $b_0$; there exist two vertex-disjoint paths $P_{11}$ and $P_{12}$ cover all vertices of $B^1-F$, where $P_{11}$ connects $u_1$ to $t_1$ and $P_{12}$ connects $v_1$ to $t_2$; there exist two vertex-disjoint paths $P_{21}$ and $P_{22}$ cover all vertices of $B^2-F$, where $P_{21}$ connects $u_2$ to $a_2$ and $P_{22}$ connects $v_2$ to $b_2$; there exist two vertex-disjoint paths $P_{31}$ and $P_{32}$ cover all vertices of $B^3-F$, where $P_{31}$ connects $u_3$ to $a_3$ and $P_{32}$ connects $v_3$ to $b_3$. Hence, $\langle s_1,P_{01},a_0,u_3,P_{31},a_3,u_2,P_{21},a_2,u_1,P_{11},t_1\rangle$ and $\langle s_2,P_{02},b_0,v_3,P_{32},b_3,v_2,P_{22},b_2,v_1,P_{12},t_1\rangle$ are two vertex-disjoint paths required.


\noindent{\bf Case 2.1.6.} $s_2\in V(B^0)$ and $t_1,t_2\in V(B^2)$. By Lemma \ref{exist-HP-abcd}, there exist four vertices $a,c\in V_0$ and $b,d\in V_1$ of $B^3$ such that:

(1) $a=p(c)$, $b=p(d)$ and $a,b,c$ and $d$ form a 4-cycle in $B^3$;

(2) there exists an $(n-1)$-dimension neighbor $a_{0}$ of $a$ and $c$ such that $a_{0}a,a_0c\not\in F$;

(3) there exist two $(n-1)$-dimension neighbors $u_{2}$ and $v_{2}$ of $b$ and $d$ such that $u_2b,u_2d,v_2b\not\in F$ and $cd\not\in F$;

(4) there exists a neighbor $u$ of $b$ and $d$ in $B^3$ such that $ub_0\not\in F$ is an $(n-1)$-dimension edge;

(5) there exists a longest path $P_3$ from $u$ to $a$ covering all vertices of $B^3-F$ but $b,c$ and $d$.

\begin{figure}
\begin{minipage}[t]{0.5\linewidth}
\centering
\includegraphics[width=60mm]{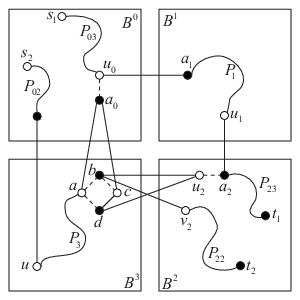}
\caption{Illustration of Cases 2.1.6.} \label{g_2_1_6}
\end{minipage}
\begin{minipage}[t]{0.5\linewidth}
\centering
\includegraphics[width=60mm]{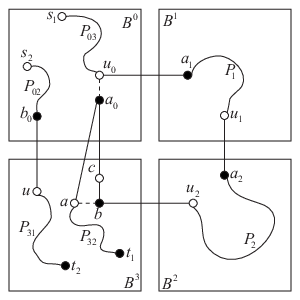}
\caption{Illustration of Case 2.1.7.} \label{g_2_1_7}
\end{minipage}
\end{figure}

It is obvious that $a_0\neq p(b_0)$. By induction hypothesis, there exist two vertex-disjoint paths $P_{01}$ and $P_{02}$ cover all vertices of $B^0-F$, where $P_{01}$ connects $s_1$ to $a_0$ and $P_{02}$ connects $s_2$ to $b_0$; there exist two vertex-disjoint paths $P_{21}$ and $P_{22}$ cover all vertices of $B^2-F$, where $P_{21}$ connects $u_2$ to $t_1$ and $P_{22}$ connects $v_2$ to $t_2$. Let $u_0$ (resp. $a_2$) be the neighbor of $a_0$ (resp. $u_2$) on $P_{01}$ (resp. $P_{21}$). For convenience, we denote $P_{01}-a_0$ by $P_{03}$, that is, $P_{03}$ is a path from $s_1$ to $u_0$. Similarly, we denote $P_{21}-a_2$ by $P_{23}$, that is, $P_{23}$ is a path from $a_2$ to $t_1$. If $|e_{n-1}(u_0)\cap F|=2$ or $|e_{n-1}(a_2)\cap F|=2$, say $|e_{n-1}(u_0)\cap F|=2$, let $u_0'\in V(B^0)$ such that $u_0'=p(u_0)$. Moreover, if $u_0'a_0\not\in F$, we can replace $u_0$ by $u_0'$ on $P_{03}$. Otherwise we have at least three fault edges incident to $u_0$ and $u_0'$. Since there are $2n-2$ common neighbors of $u_0$ and $u_0'$ in $B^0$, fault edges incident to $u_0$ and $u_0'$ may affect $2n-2$ vertices as the choice of $a_0$. Since $3\times((4^{n-1}-2)/2)/(2n-2)> 2n-3$ whenever $n\geq3$, we can always choose such $u_0\in V(B^0)$ and $a_2\in V(B^2)$ that there exist two fault-free $(n-1)$-dimension edges $u_0a_1\in E_{0,1}$ and $u_1a_2\in E_{1,2}$. Then there exists a fault-free Hamiltonian path $P_{1}$ of $B^1-F$ from $a_1$ to $u_1$. Hence, $\langle s_1,P_{03},u_0,a_1,P_{1},u_1,a_2,P_{23},t_1\rangle$ and $\langle s_2,P_{02},b_0,u,P_{3},a,a_0,c,d,u_2,b,v_2,P_{22},t_2\rangle$ are two vertex-disjoint paths required (see Fig. \ref{g_2_1_6}).

\noindent{\bf Case 2.1.7.} $s_2\in V(B^0)$ and $t_1,t_2\in V(B^3)$. By Lemma \ref{exist-HP-t_1_abc}, there exist four vertices $u,a,c\in V_0$ and $b\in V_1$ of $B^3$ with $a=p(c)$ such that:

(1) there exists an $(n-1)$-dimension neighbor $a_{0}$ of $a$ and $c$ such that $a_{0}a,a_{0}c\not\in F$ and there exists an $(n-1)$-dimension neighbor $u_{2}$ of $b$ such that $u_2b\not\in F$, where $b$ ($b\neq t_1,t_2$) is a common neighbor of $a$ and $c$;

(2) there exists an $(n-1)$-dimension neighbor $b_0$ of $u$ such that $ub_{0}\not\in F$;

(3) there exist two vertex-disjoint paths $P_{31}$ and $Q$ of $B^3-F$ cover all vertices of it, where $P_{31}$ connects $u$ to $t_2$, and $Q$ connects $c$ to $t_1$ and $\langle c,b,a\rangle$ is a subpath of $Q$.

Deleting $ab$ from $Q$ will generate two vertex-disjoint paths $bc$ and $P_{32}$, where $P_{32}$ connects $a$ to $t_1$. By induction hypothesis, there exist two vertex-disjoint paths $P_{01}$ and $P_{02}$ cover all vertices of $B^0-F$, where $P_{01}$ connects $s_1$ to $a_0$ and $P_{02}$ connects $s_2$ to $b_0$. Similar to the proof of Case 2.1.6, let $u_0$ be the neighbor of $a_0$ on $P_{01}$ such that $u_0a_1\in E_{0,1}$ is a fault-free edge. For convenience, we denote $P_{01}-a_0$ by $P_{03}$, that is, $P_{03}$ is a path from $s_1$ to $u_0$. By Lemma \ref{crossing-edges}, there must exist a fault-free edge $u_1a_2\in E_{1,2}$. Additionally, there exist a fault-free Hamiltonian path $P_1$ of $B^1-F$ from $a_1$ to $u_1$, and a fault-free Hamiltonian path $P_2$ of $B^2-F$ from $a_2$ to $u_2$. Hence, $\langle s_1,P_{03},u_0,a_1,P_{1},u_1,a_2,P_{2},u_2,b,c,a_0,a,P_{32},t_1\rangle$ and $\langle s_2,P_{02},b_0,u,P_{31},t_2\rangle$ are two vertex-disjoint paths required (see Fig. \ref{g_2_1_7}).


\noindent{\bf Case 2.2.} For all $j\in \{0,1,2\}$, $|V(B^j)\cap \{s_2,t_1,t_2\}|\leq1$.

\noindent{\bf Case 2.2.1.} $t_1\in V(B^0)$.

\noindent{\bf Case 2.2.1.1.} $t_2\in V(B^1)$ and $s_2\in V(B^2)$. There exist fault-free edges $u_0a_1\in E_{0,1}$ such that $u_0\neq s_1$, $u_1a_2,v_1b_2\in E_{1,2}$ such that $u_1\neq v_1$ and $a_2\neq b_2$, $u_2a_3\in E_{2,3}$ such that $u_2\neq s_2$, and $u_3a_0\in E_{3,0}$ such that $a_0\neq t_1$. By induction hypothesis, there exist two vertex-disjoint paths $P_{01}$ and $P_{02}$ cover all vertices of $B^0-F$, where $P_{01}$ connects $s_1$ to $a_0$ and $P_{02}$ connects $u_0$ to $t_1$; there exist two vertex-disjoint paths $P_{11}$ and $P_{12}$ cover all vertices of $B^1-F$, where $P_{11}$ connects $u_1$ to $a_1$ and $P_{12}$ connects $v_1$ to $t_2$; there exist two vertex-disjoint paths $P_{21}$ and $P_{22}$ cover all vertices of $B^2-F$, where $P_{21}$ connects $u_2$ to $a_2$ and $P_{22}$ connects $s_2$ to $b_2$. By Lemma \ref{HP-fault}, there exists a fault-free Hamiltonian path $P_3$ of $B^3-F$ from $u_3$ to $a_3$. Hence, $\langle s_1,P_{01},a_0,u_3,P_{3},a_3,u_2,P_{21},a_2,u_1,P_{11},a_1,u_0,P_{02},t_1\rangle$ and $\langle s_2,P_{22},b_2,v_1,P_{12},t_2\rangle$ are two vertex-disjoint paths required.

\noindent{\bf Case 2.2.1.2.} $t_2\in V(B^1)$ and $s_2\in V(B^3)$. There exist two non-faulty edges $u_1a_2\in E_{1,2}$ and $u_2a_3\in E_{2,3}$. By Lemma \ref{HP-fault}, there exist a fault-free Hamiltonian path $P_0$ of $B^0-F$ from $s_1$ to $t_1$, a fault-free Hamiltonian path $P_1$ of $B^1-F$ from $u_1$ to $t_2$, a fault-free Hamiltonian path $P_2$ of $B^2-F$ from $u_2$ to $a_2$, and a fault-free Hamiltonian path $P_3$ of $B^3-F$ from $s_2$ to $a_3$. Hence, $\langle s_1,P_{1},t_1\rangle$ and $\langle s_2,P_{3},a_3,u_2,P_{2},a_2,u_1,P_{1},t_2\rangle$ are two vertex-disjoint paths required.

\noindent{\bf Case 2.2.1.3.} $t_2\in V(B^2)$ and $s_2\in V(B^3)$. There exist fault-free edges $u_0a_1\in E_{0,1}$ such that $u_0\neq s_1$, $u_1a_2\in E_{1,2}$ such that $a_2\neq t_2$, $u_2a_3,v_2b_3\in E_{2,3}$ such that $a_3\neq b_3$ and $u_2\neq v_2$, and $u_3a_0\in E_{2,3}$ such that $u_3\neq s_2$ and $a_0\neq t_1$. By induction hypothesis, there exist two vertex-disjoint paths $P_{01}$ and $P_{02}$ cover all vertices of $B^0-F$, where $P_{01}$ connects $s_1$ to $a_0$ and $P_{02}$ connects $u_0$ to $t_1$; there exist two vertex-disjoint paths $P_{21}$ and $P_{22}$ cover all vertices of $B^2-F$, where $P_{21}$ connects $u_2$ to $a_2$ and $P_{22}$ connects $v_2$ to $t_2$; there exist two vertex-disjoint paths $P_{31}$ and $P_{32}$ cover all vertices of $B^3-F$, where $P_{31}$ connects $u_3$ to $a_3$ and $P_{32}$ connects $s_2$ to $b_3$. By Lemma \ref{HP-fault}, there exists a fault-free Hamiltonian path $P_1$ of $B^1-F$ from $u_1$ to $a_1$. Hence, $\langle s_1,P_{01},a_0,u_3,P_{31},a_3,u_2,P_{21},a_2,u_1,P_{1},a_1,u_0,P_{02},t_1\rangle$ and $\langle s_2,P_{32},b_3,v_2,P_{22},t_2\rangle$ are two vertex-disjoint paths required.

\noindent{\bf Case 2.2.1.4.} $t_2\in V(B^2)$ and $s_2\in V(B^1)$. There exist fault-free edges $u_0a_1\in E_{0,1}$ such that $u_0\neq s_1$, $v_1b_2\in E_{1,2}$ such that $v_1\neq s_2$ and $b_2\neq t_2$, $u_2a_3,v_2b_3\in E_{2,3}$ such that $a_3\neq b_3$ and $u_2\neq v_2$, and $u_3a_0\in E_{3,0}$ such that $a_0\neq t_1$. By induction hypothesis, there exist two vertex-disjoint paths $P_{01}$ and $P_{02}$ cover all vertices of $B^0-F$, where $P_{01}$ connects $u_0$ to $a_0$ and $P_{02}$ connects $s_1$ to $t_1$; there exist two vertex-disjoint paths $P_{21}$ and $P_{22}$ cover all vertices of $B^2-F$, where $P_{21}$ connects $v_2$ to $b_2$ and $P_{22}$ connects $u_2$ to $t_2$. Moreover, there must exist an edge $v_0b_0$ in $P_{01}$ or $P_{02}$, say $P_{02}$ such that there exist two fault-free $(n-1)$-dimension edges $v_0b_1\in E_{0,1}$ and $v_3b_0\in E_{3,0}$, where $v_3\neq u_3$ and $b_1\neq a_1$. Deleting $v_0b_0$ from $P_{02}$ will generate two vertex-disjoint paths $P_{03}$ and $P_{04}$, where $P_{03}$ connects $s_1$ to $b_0$ and $P_{04}$ connects $v_0$ to $t_1$. Analogously, there exist two vertex-disjoint paths $P_{11}$ and $P_{12}$ cover all vertices of $B^1-F$, where $P_{11}$ connects $v_1$ to $b_1$ and $P_{12}$ connects $s_2$ to $a_1$; there exist two vertex-disjoint paths $P_{31}$ and $P_{32}$ cover all vertices of $B^3-F$, where $P_{31}$ connects $v_3$ to $b_3$ and $P_{32}$ connects $u_3$ to $a_3$. Hence, $\langle s_1,P_{03},b_0,v_3,P_{31},b_3,v_2,P_{21},b_2,v_1,P_{11},b_1,v_0,P_{04},t_1\rangle$ and $\langle s_2,P_{12},a_1,u_0,P_{01},a_0,u_3,P_{32},a_3,u_2,P_{22},t_2\rangle$ are two vertex-disjoint paths required.

\noindent{\bf Case 2.2.1.5.} $t_2\in V(B^3)$ and $s_2\in V(B^1)$. The proof is quite analogous to that of Case 2.2.1.4, we omit it.

\noindent{\bf Case 2.2.1.6.} $t_2\in V(B^3)$ and $s_2\in V(B^2)$. The proof is quite analogous to that of Case 2.2.1.4, we omit it.

\noindent{\bf Case 2.2.2.} $t_2\in V(B^0)$.

\noindent{\bf Case 2.2.2.1.} $t_1\in V(B^1)$ and $s_2\in V(B^2)$. There exist fault-free edges $v_0b_1\in E_{0,1}$ such that $v_0\neq s_1$ and $b_1\neq t_1$, $u_1a_2,v_1b_2\in E_{1,2}$ such that $u_1\neq v_1$ and $a_2\neq b_2$, $u_2a_3\in E_{2,3}$ such that $u_2\neq s_2$, and $u_3a_0\in E_{3,0}$ such that $a_0\neq t_2$. By induction hypothesis, there exist two vertex-disjoint paths $P_{01}$ and $P_{02}$ cover all vertices of $B^0-F$, where $P_{01}$ connects $s_1$ to $a_0$ and $P_{02}$ connects $v_0$ to $t_2$; there exist two vertex-disjoint paths $P_{11}$ and $P_{12}$ cover all vertices of $B^1-F$, where $P_{11}$ connects $u_1$ to $t_1$ and $P_{12}$ connects $v_1$ to $b_1$; there exist two vertex-disjoint paths $P_{21}$ and $P_{22}$ cover all vertices of $B^2-F$, where $P_{21}$ connects $u_2$ to $a_2$ and $P_{22}$ connects $s_2$ to $b_2$. By Lemma \ref{HP-fault}, there exists a Hamiltonian path $P_{3}$ of $B^3-F$ from $u_3$ to $a_3$. Hence, $\langle s_1,P_{01},a_0,u_3,P_{3},a_3,u_2,P_{21},a_2,u_1,P_{11},t_1\rangle$ and $\langle s_2,P_{22},b_2,v_1,P_{12},b_1,v_0,P_{02},t_2\rangle$ are two vertex-disjoint paths required.

\noindent{\bf Case 2.2.2.2.} $t_1\in V(B^1)$ and $s_2\in V(B^3)$. There exist fault-free edges $v_0b_1\in E_{0,1}$ such that $v_0\neq s_1$ and $b_1\neq t_1$, $u_1a_2,v_1b_2\in E_{1,2}$ such that $u_1\neq v_1$ and $a_2\neq b_2$, $u_2a_3,v_2b_3\in E_{2,3}$ such that $u_2\neq v_2$ and $a_3\neq b_3$, and $u_3a_0\in E_{3,0}$ such that $a_0\neq t_2$ and $u_3\neq s_2$. By induction hypothesis, there exist two vertex-disjoint paths $P_{01}$ and $P_{02}$ cover all vertices of $B^0-F$, where $P_{01}$ connects $s_1$ to $a_0$ and $P_{02}$ connects $v_0$ to $t_2$; there exist two vertex-disjoint paths $P_{11}$ and $P_{12}$ cover all vertices of $B^1-F$, where $P_{11}$ connects $u_1$ to $t_1$ and $P_{12}$ connects $v_1$ to $b_1$; there exist two vertex-disjoint paths $P_{21}$ and $P_{22}$ cover all vertices of $B^2-F$, where $P_{21}$ connects $u_2$ to $a_2$ and $P_{22}$ connects $v_2$ to $b_2$; there exist two vertex-disjoint paths $P_{31}$ and $P_{32}$ cover all vertices of $B^3-F$, where $P_{31}$ connects $u_3$ to $a_3$ and $P_{32}$ connects $s_2$ to $b_3$. Hence, $\langle s_1,P_{01},a_0,u_3,P_{31},a_3,u_2,P_{21},a_2,u_1,P_{11},t_1\rangle$ and $\langle s_2,P_{32},b_3,v_2,P_{22},b_2,v_1,P_{12},b_1,v_0,P_{02},t_2\rangle$ are two vertex-disjoint paths required.

\noindent{\bf Case 2.2.2.3.} $t_1\in V(B^2)$ and $s_2\in V(B^1)$. There exist fault-free edges $v_0b_1\in E_{0,1}$ such that $v_0\neq s_1$, $u_2a_3\in E_{2,3}$, and $u_3a_0\in E_{3,0}$ such that $a_0\neq t_2$. By induction hypothesis, there exist two vertex-disjoint paths $P_{01}$ and $P_{02}$ cover all vertices of $B^0-F$, where $P_{01}$ connects $s_1$ to $a_0$ and $P_{02}$ connects $v_0$ to $t_2$. By Lemma \ref{HP-fault}, there exist a fault-free Hamiltonian path $P^1$ of $B^1-F$ from $s_2$ to $b_1$, a fault-free Hamiltonian path $P^2$ of $B^2-F$ from $u_2$ to $t_1$, and a fault-free Hamiltonian path $P^3$ of $B^3-F$ from $u_3$ to $a_3$. Hence, $\langle s_1,P_{01},a_0,u_3,P_{3},a_3,u_2,P_{2},t_1\rangle$ and $\langle s_2,P_{1},b_1,v_0,P_{02},t_2\rangle$ are two vertex-disjoint paths required.

\noindent{\bf Case 2.2.2.4.} $t_1\in V(B^2)$ and $s_2\in V(B^3)$. The proof is quite analogous to that of Case 2.2.2.1, we omit it.

\noindent{\bf Case 2.2.2.5.} $t_1\in V(B^3)$ and $s_2\in V(B^1)$. There exist fault-free edges $v_0b_1\in E_{0,1}$ such that $v_0\neq s_1$, $u_1a_2\in E_{1,2}$ such that $u_1\neq s_2$, $u_2a_3\in E_{2,3}$ such that $a_3\neq t_1$, and $v_3a_0\in E_{3,0}$ such that $a_0\neq t_2$. By induction hypothesis, there exist two vertex-disjoint paths $P_{01}$ and $P_{02}$ cover all vertices of $B^0-F$, where $P_{01}$ connects $v_0$ to $t_2$ and $P_{02}$ connects $s_1$ to $a_0$. Moreover, there must exist an edge $u_0b_0$ on $P_{01}$ or $P_{02}$, say $P_{02}$ such that there exist two fault-free $(n-1)$-dimension edges $u_0a_1\in E_{0,1}$ and $u_3b_0\in E_{3,0}$, where $u_3\neq v_3$ and $a_1\neq b_1$. Deleting $u_0b_0$ from $P_{02}$ will generate two vertex-disjoint paths $P_{03}$ and $P_{04}$, where $P_{03}$ connects $s_1$ to $b_0$ and $P_{04}$ connects $u_0$ to $a_0$. Analogously, there exist two vertex-disjoint paths $P_{11}$ and $P_{12}$ cover all vertices of $B^1-F$, where $P_{11}$ connects $u_1$ to $a_1$ and $P_{12}$ connects $s_2$ to $b_1$; there exist two vertex-disjoint paths $P_{31}$ and $P_{32}$ cover all vertices of $B^3-F$, where $P_{31}$ connects $v_3$ to $t_1$ and $P_{32}$ connects $u_3$ to $a_3$. By Lemma \ref{HP-fault}, there exist a fault-free Hamiltonian path $P_2$ of $B^2-F$ from $u_2$ to $a_2$. Hence, $\langle s_1,P_{03},b_0,u_3,P_{32},a_3,u_2,P_{2},a_2,u_1,P_{11},a_1,u_0,P_{04},a_0,v_3,P_{31},t_1\rangle$ and $\langle s_2,P_{12},b_1,v_0,P_{01},t_2\rangle$ are two vertex-disjoint paths required.

\noindent{\bf Case 2.2.2.6.} $t_1\in V(B^3)$ and $s_2\in V(B^2)$. There exist fault-free edges $v_0b_1\in E_{0,1}$ such that $v_0\neq s_1$, $v_1b_2\in E_{1,2}$, and $u_3a_0\in E_{3,0}$ such that $a_0\neq t_2$. By induction hypothesis, there exist two vertex-disjoint paths $P_{01}$ and $P_{02}$ cover all vertices of $B^0-F$, where $P_{01}$ connects $s_1$ to $a_0$ and $P_{02}$ connects $v_0$ to $t_2$. By Lemma \ref{HP-fault}, there exist a fault-free Hamiltonian path $P_1$ of $B^1-F$ from $v_1$ to $b_1$, a fault-free Hamiltonian path $P_2$ of $B^2-F$ from $s_2$ to $b_2$, and a fault-free Hamiltonian path $P_3$ of $B^3-F$ from $u_3$ to $t_1$. Hence, $\langle s_1,P_{01},a_0,u_3,P_{3},t_1\rangle$ and $\langle s_2,P_{2},b_2,v_1,P_{1},b_1,v_0,P_{02},t_2\rangle$ are two vertex-disjoint paths required.

\noindent{\bf Case 2.2.3.} $s_2\in V(B^0)$.

\noindent{\bf Case 2.2.3.1.} $t_1\in V(B^1)$ and $t_2\in V(B^2)$. There exist fault-free edges $u_1a_2\in E_{1,2}$ such that $a_2\neq t_2$, $u_2a_3,v_2b_3\in E_{2,3}$ such that $u_2\neq v_2$ and $a_3\neq b_3$, and $u_3a_0,v_3b_0\in E_{3,0}$ such that $u_3\neq v_3$ and $a_0\neq b_0$. By induction hypothesis, there exist two vertex-disjoint paths $P_{01}$ and $P_{02}$ cover all vertices of $B^0-F$, where $P_{01}$ connects $s_1$ to $a_0$ and $P_{02}$ connects $s_2$ to $b_0$; there exist two vertex-disjoint paths $P_{21}$ and $P_{22}$ cover all vertices of $B^2-F$, where $P_{21}$ connects $u_2$ to $a_2$ and $P_{22}$ connects $v_2$ to $t_2$; there exist two vertex-disjoint paths $P_{31}$ and $P_{32}$ cover all vertices of $B^3-F$, where $P_{31}$ connects $u_3$ to $a_3$ and $P_{32}$ connects $v_3$ to $b_3$. By Lemma \ref{HP-fault}, there exists a fault-free Hamiltonian path $P_{1}$ of $B^1-F$ from $u_1$ to $t_1$. Hence, $\langle s_1,P_{01},a_0,u_3,P_{31},a_3,u_2,P_{21},a_2,u_1,P_{1},t_1\rangle$ and $\langle s_2,P_{02},b_0,v_3,P_{32},b_3,v_2,P_{22},t_2\rangle$ are two vertex-disjoint paths required.

\noindent{\bf Case 2.2.3.2.} $t_1\in V(B^1)$ and $t_2\in V(B^3)$. The proof is quite analogous to that of Case 2.2.3.1, we omit it.

\noindent{\bf Case 2.2.3.3.} $t_1\in V(B^2)$ and $t_2\in V(B^3)$. By Lemma \ref{exist-HP-t_1_abc}, there exist four vertices $u,a,c\in V_1$ and $b\in V_0$ of $B^0-F$ with $a=p(c)$ such that:

(1) there exists an $(n-1)$-dimension neighbor $u_{3}$ of $a$ and $c$ such that $u_{3}a,u_{3}c\not\in F$ and there exists an $(n-1)$-dimension neighbor $a_{1}$ of $b$ such that $a_1b\not\in F$, where $b$ ($b\neq s_1,s_2$) is a common neighbor of $a$ and $c$;

(2) there exists an $(n-1)$-dimension neighbor $v_3$ of $u$ such that $uv_{3}\not\in F$;

(3) there exist two vertex-disjoint paths $P_{01}$ and $Q$ of $B^0-F$ cover all vertices of it, where $P_{01}$ connects $s_2$ to $u$, and $Q$ connects $s_1$ to $c$ and $\langle c,b,a\rangle$ is a subpath of $Q$.

Deleting $ab$ from $Q$ will generate two vertex-disjoint paths $P_{02}$ and $bc$, where $P_{02}$ connects $s_1$ to $a$ and $bc$ is an edge. Let $a_3\in V_1$ be a vertex in $B^3$ such that $a_3\neq t_2$ and $u_2a_3\in E_{2,3}$ is a fault-free edge. In addition, there exist two vertex-disjoint paths $P_{31}$ and $P_{32}$ cover all vertices of $B^3-F$, where $P_{31}$ connects $v_3$ to $t_2$ and $P_{32}$ connects $u_3$ to $a_3$. Similar to the proof of Case 2.1.6, let $b_3$ be the neighbor of $u_3$ on $P_{32}$ such that $v_2b_3\in E_{2,3}$ is a fault-free edge. For convenience, we denote $P_{32}-u_3$ by $P_{33}$, that is, $P_{33}$ is a path from $b_3$ to $a_3$. By Lemma \ref{crossing-edges}, there must exist an fault-free edge $u_1a_2\in E_{1,2}$ such that $a_2\neq t_1$. Thus, there exist two vertex-disjoint paths $P_{21}$ and $P_{22}$ cover all vertices of $B^2-F$, where $P_{21}$ connects $u_2$ to $t_1$ and $P_{22}$ connects $a_2$ to $v_2$. Additionally, there exist a fault-free Hamiltonian path $P_1$ of $B^1-F$ from $a_1$ to $u_1$. Hence, $\langle s_1,P_{02},a,u_3,c,b,a_1,P_{1},u_1,a_2,P_{22},v_2,b_3,P_{33},a_3,u_2,P_{21},t_1\rangle$ and $\langle s_2,P_{01},u,v_3,P_{31},t_2\rangle$ are two vertex-disjoint paths required.



\noindent{\bf Case 3.} $|V(B^0)\cap \{s_2,t_1,t_2\}|=2$.

\noindent{\bf Case 3.1.} $t_1,t_2\in V(B^0)$ and $s_2\in V(B^1)$. There exist a fault-free edge $v_0b_1\in E_{0,1}$ such that $v_0\neq s_1$. By induction hypothesis, there exist two vertex-disjoint paths $P_{01}$ and $P_{02}$ cover all vertices of $B^0-F$, where $P_{01}$ connects $v_0$ to $t_2$ and $P_{02}$ connects $s_1$ to $t_1$. Moreover, there must exist an edge $u_0a_0$ on $P_{01}$ or $P_{02}$, say $P_{02}$ such that there exist two fault-free $(n-1)$-dimension edges $u_0a_1\in E_{0,1}$ and $u_3a_0\in E_{3,0}$, where $a_1\neq b_1$. Deleting $u_0a_0$ from $P_{02}$ will generate two vertex-disjoint paths $P_{03}$ and $P_{04}$, where $P_{03}$ connects $s_1$ to $a_0$ and $P_{04}$ connects $u_0$ to $t_1$. In addition, there exist two fault-free edges $u_1a_2\in E_{1,2}$ and $u_2a_3\in E_{2,3}$, where $u_1\neq s_2$ . Analogously, there exist two vertex-disjoint paths $P_{11}$ and $P_{12}$ cover all vertices of $B^1-F$, where $P_{11}$ connects $u_1$ to $a_1$ and $P_{12}$ connects $s_2$ to $b_1$. Moreover, by Lemma \ref{HP-fault}, there exist a fault-free Hamiltonian path $P_2$ of $B^2-F$ from $u_2$ to $a_2$ and a fault-free Hamiltonian path $P_3$ of $B^3-F$ from $u_3$ to $a_3$. Hence, $\langle s_1,P_{03},a_0,u_3,P_{3},a_3,u_2,P_{2},a_2,u_1,P_{11},a_1,u_0,P_{04},t_1\rangle$ and $\langle s_2,P_{12},b_1,v_0,P_{01},t_2\rangle$ are two vertex-disjoint paths required.

\noindent{\bf Case 3.2.} $t_1,t_2\in V(B^0)$ and $s_2\in V(B^2)$. There exist a fault-free edge $v_0b_1\in E_{0,1}$ such that $v_0\neq s_1$. By induction hypothesis, there exist two vertex-disjoint paths $P_{01}$ and $P_{02}$ cover all vertices of $B^0-F$, where $P_{01}$ connects $v_0$ to $t_2$ and $P_{02}$ connects $s_1$ to $t_1$. Moreover, there must exist an edge $u_0a_0$ on $P_{01}$ or $P_{02}$, say $P_{02}$ such that there exist two fault-free $(n-1)$-dimension edges $u_0a_1\in E_{0,1}$ and $u_3a_0\in E_{3,0}$, where $a_1\neq b_1$. Deleting $u_0a_0$ from $P_{02}$ will generate two vertex-disjoint paths $P_{03}$ and $P_{04}$, where $P_{03}$ connects $s_1$ to $a_0$ and $P_{04}$ connects $u_0$ to $t_1$. In addition, there exist fault-free edges $u_1a_2,v_1b_2\in E_{1,2}$ such that $u_1\neq v_1$ and $a_2\neq b_2$, and $u_2a_3\in E_{2,3}$. Analogously, there exist two vertex-disjoint paths $P_{11}$ and $P_{12}$ cover all vertices of $B^1-F$, where $P_{11}$ connects $u_1$ to $a_1$ and $P_{12}$ connects $v_1$ to $b_1$; there exist two vertex-disjoint paths $P_{21}$ and $P_{22}$ cover all vertices of $B^2-F$, where $P_{21}$ connects $u_2$ to $a_2$ and $P_{22}$ connects $s_2$ to $b_2$. Moreover, by Lemma \ref{HP-fault}, there exist a fault-free Hamiltonian path $P_3$ of $B^3-F$ from $u_3$ to $a_3$. Hence, $\langle s_1,P_{03},a_0,u_3,P_{3},a_3,u_2,P_{21},a_2,u_1,P_{11},a_1,u_0,P_{04},t_1\rangle$ and $\langle s_2,P_{22},b_2,v_1,P_{12},b_1,v_0,P_{01},t_2\rangle$ are two vertex-disjoint paths required.

\noindent{\bf Case 3.3.} $t_1,t_2\in V(B^0)$ and $s_2\in V(B^3)$. There exist fault-free edges $v_0b_1\in E_{0,1}$ such that $v_0\neq s_1$, $v_1b_2\in E_{1,2}$, and $v_2b_3\in E_{2,3}$. By induction hypothesis, there exist two vertex-disjoint paths $P_{01}$ and $P_{02}$ cover all vertices of $B^0-F$, where $P_{01}$ connects $s_1$ to $t_1$ and $P_{02}$ connects $v_0$ to $t_2$. By Lemma \ref{HP-fault}, there exist a fault-free Hamiltonian path $P_1$ of $B^1-F$ from $v_1$ to $b_1$, a fault-free Hamiltonian path $P_2$ of $B^2-F$ from $v_2$ to $b_2$, and a fault-free Hamiltonian path $P_3$ of $B^3-F$ from $s_2$ to $b_3$. Hence, $\langle s_1,P_{01},t_1\rangle$ and $\langle s_2,P_{3},b_3,v_2,P_{2},b_2,v_1,P_{1},b_1,v_0,P_{02},t_2\rangle$ are two vertex-disjoint paths required.

\noindent{\bf Case 4.} $s_2,t_1,t_2\in V(B^0)$. By induction hypothesis, there exist two vertex-disjoint paths $P_{01}$ and $P_{02}$ cover all vertices of $B^0-F$, where $P_{01}$ connects $s_1$ to $t_1$ and $P_{02}$ connects $s_2$ to $t_2$. Since $l(P_{01})+l(P_{02})=4^{n-1}-2$ and any vertex has two $(n-1)$-dimension neighbors, there must exist an edge $u_0a_0$ on $P_{01}$ or $P_{02}$, say $P_{01}$ such that there exist two non-faulty edges $u_0a_1\in E_{0,1}$ and $u_3a_0\in E_{3,0}$. Thus, deleting $u_0a_0$ from $P_{01}$ will generate two vertex-disjoint paths $P_{03}$ and $P_{04}$, where where $P_{03}$ connects $s_1$ to $a_0$ and $P_{04}$ connects $u_0$ to $t_1$. In addition, there exist non-faulty edges $u_1a_2\in E_{1,2}$ and $u_2a_3\in E_{2,3}$. By Lemma \ref{HP-fault}, there exist a fault-free Hamiltonian path $P_1$ of $B^1-F$ from $u_1$ to $a_1$, a fault-free Hamiltonian path $P_2$ of $B^2-F$ from $u_2$ to $a_2$, and a fault-free Hamiltonian path $P_3$ of $B^3-F$ from $u_3$ to $a_3$. Hence, $\langle s_1,P_{03},a_0,u_3,P_{3},a_3,u_2,P_{2},a_2,u_1,P_{1},a_1,u_0,P_{04},t_1\rangle$ and $\langle s_2,P_{02},t_2\rangle$ are two vertex-disjoint paths required.

Thus, this completes the proof.\qed

\section{Conclusions}
In this paper, we consider paired two-disjoint path cover of the balanced hypercube with some faulty edges. We use induction to prove that the balanced hypercube $BH_n$, $n\geq2$, is paired two-disjoint path coverable when at most $2n-3$ edge faults occur.

\vskip 0.1 in

\begin{figure}[h]
\centering
\includegraphics[width=46mm]{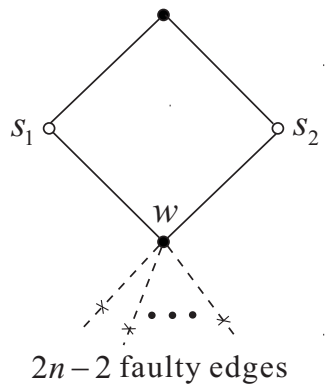}
\caption{$BH_n$ has no paired two-disjoint path cover with $2n-2$ faulty edges.} \label{g_conclusion}
\end{figure}

Let $s_1,s_2\in V_0$ and $t_1,t_2\in V_1$ be four vertices in $BH_n$. There exists a balanced hypercube $BH_n$ with $2n-2$ edge faults containing no vertex-disjoint paths $P_i$, $i=1,2$, that cover all vertices of it, where $P_i$ connects $s_i$ to $t_i$ and $V(P_1)\cup V(P_2)=V(BH_n)$. For example, let $s_1$ and $s_2$ be two vertices differing only from the inner index and let $w$ be any common neighbor of $s_1$ and $s_2$. One can consider that the $2n-2$ edges incident to $w$ (except $s_1w$ and $s_2w$) are all faulty (see Fig. \ref{g_conclusion}). Obviously, $w$ has exactly two fault-free edges incident to it. Therefore, it is impossible to have vertex-disjoint $s_1,t_1$-path and $s_2,t_2$-path that cover all vertices of $BH_n$. Hence, our results are optimal.

\noindent{\bf\large Acknowledgments}

I am grateful to Simon R. Blackburn for fruitful discussions during my visit to Royal Holloway, University of London.

\end{document}